\documentclass[11pt]{amsart}
\usepackage{geometry}                
\usepackage{graphicx}
\usepackage{amssymb}
\usepackage{epstopdf}
\newtheorem{lemma}{Lemma}[subsection]
\newtheorem{propos}[lemma]{Proposition}

\newtheorem{theorem}[lemma]{Theorem}
\newtheorem{cor}[lemma]{Corollary}
\newtheorem{corol}[lemma]{Corollary}
\newtheorem{defin}[lemma]{Definition}
\newtheorem{remark}[lemma]{Remark}

\newcommand{\CM}{\hbox{{$\mathcal M$}}}

\newcommand{\CE}{\hbox{{$\mathcal E$}}}

\newcommand{\C}{\mathbb{C}}
\newcommand{\R}{\mathbb{R}}

\newcommand{\Z}{\mathbb{Z}}

\newcommand{\Hom}{\mathrm{Hom}}
\renewcommand{\proof}{{\noindent {\bfseries  Proof:}\quad }}

\newcommand{\ev}{\mathrm{ev}}
\newcommand{\coev}{\mathrm{coev}}

\newcommand{\extd}{\mathrm{d}}

\newcommand{\isom}{{\cong}}
\newcommand{\eps}{{\epsilon}}
\newcommand{\tens}{\mathop{\otimes}}

\newcommand{\ra}{{\triangleleft}}

\newcommand{\Ad}{\mathrm{ Ad}}

\newcommand{\id}{\mathrm{id}}

\newcommand{\<}{\langle}
\renewcommand{\>}{\rangle}

\renewcommand{\o}{{}_{\scriptscriptstyle(1)}}
\renewcommand{\t}{{}_{\scriptscriptstyle(2)}}
\newcommand{\thr}{{}_{\scriptscriptstyle(3)}}

\title{$*$-Compatible Connections in Noncommutative Riemannian Geometry}
\author{E.J. Beggs  \& S. Majid}
\address{EJB: Department of Mathematics,  Swansea University, Singleton Parc, SA2 8PP, SM: School of Mathematical Sciences, Queen Mary University of London, Mile End Rd, London E1, 4NS}

\begin{document}

\abstract{We develop the formalism for noncommutative differential geometry and Riemmannian geometry  to take full account of the $*$-algebra structure on the (possibly noncommutative) coordinate ring and the  bimodule structure on the differential forms. We show that $*$-compatible bimodule connections lead to braid operators $\sigma$ in some generality (going beyond the quantum group case) and we develop their role in the exterior algebra. We study metrics in the form of Hermitian structures  on Hilbert $*$-modules and metric compatibility in both the usual and a cotorsion form. We show that the theory works well for the quantum group $\C_q[SU_2]$ with its 3D calculus, finding for each point of a 3-parameter space of covariant metrics a unique `Levi-Civita' connection deforming the classical one and characterised by zero torsion, metric-preservation and $*$-compatibility.  Allowing torsion, we find a unique connection with classical limit that is metric-preserving and $*$-compatible and for which $\sigma$ obeys the braid relations. It projects to a unique `Levi-Civita' connection on the quantum sphere.   The theory also works for finite groups and in particular  for the permutation group $S_3$ where we find somewhat similar results.}}

\maketitle

\section{Introduction}
In this paper we will consider Riemannian geometry on
some noncommutative spaces. To do this, we must reconcile the $C^*$
algebra viewpoint, where normed algebras over $\mathbb{C}$ with conjugate-linear
involutions are identified with noncommutative spaces, and the
usual point of view of classical differential geometry in terms of connections and differential forms. There is already
a well known theory of differential calculi on noncommutative algebras which includes $*$ as a differential graded $*$-algebra structure on the differential forms, but beyond this there are several different approaches to the actual Riemannian geometry. The one of Connes \cite{Connes}  based on the Dirac operator generalised as a spectral triple has, in particular, been complemented by a more constructive algebraic approach to connections, curvature, frame bundles etc., led by examples coming out of quantum group theory, see notably \cite{Ma:rief,Ma:rsph}. The latter work has typically been algebraic over a general field whereas  properly to apply it to Riemannian geometry, with its positivity requirements, requires that
we take the star structure over $\C$ fully into consideration. We provide such a formulation now,  in Sections~2 and 3. One feature of combining covariant derivatives and star operations is that
we are more or less forced to use the notion of a bimodule covariant derivative
 in which a linear connection is supplemented by some kind of `generalised braiding' $\sigma$ on the space of 1-forms. This idea had its origins in \cite{DVMic,DVMass} and \cite{Mourad}, and was later used in \cite{Madore,FioMad}. In \cite{BMHDSfinitegrps} it was shown that this idea allowed tensoringof bimodules with connections, and we refer to \cite{DVlectures} for an overview and further references. Section~2 begins with a full theory inspired in part by this earlier work, but now developed in full using our theory of bar categories.

In the new formalism, the Riemannian metric will be implemented not by a 2-covector $g$ but by the standard notion of a  Hilbert $C^*$ module (though we may omit the completeness assumption). One should think of the metric now as a sesquilinear extension of a real 2-covector, i.e. a Hermitian metric, and it is this  which turns out to be more central. We are then able to use the recently developed  language of bar categories \cite{barcats}, which makes it possible to establish conditions for
connections to be compatible with the star operations, and with Hilbert $C^*$ modules. The central issue is that $\nabla g$ naively defined by $\nabla$ extending as a derivation does not make sense when the coordinate algebra is noncommutative.  This is resolved in the Hilbert $C^*$ module approach and provides an alternative  to a weaker metric compatibility notion of `cotorsion free' previously used to address this problem\cite{Ma:rief}. We will also need and find a hermitian version of the latter weaker notion. The reader should note that for our examples we only require a fairly simple bar category;
 there are examples of rather more exotic ones in \cite{barcats}, for example relating to quantum groups at roots of unity.   The bar category language in the present paper can be thought of mostly as a book keeping device, ensuring that the  right formulae are applied in the right places. However, our approach also in principle extends to other bar categories.

The second half of the paper is concerned with showing that this general theory works well with the key `test cases' of quantum $SU_2$, quantum spheres and finite groups such as the permutation group $S_3$. We both see how previous work on these examples for specific metrics embeds into our general theory, and we gain some idea of what else the theory allows in the moduli space of connections subject to our various requirements. Section 4 provides some general remarks for the restriction of our theory to the case of Hopf algebras (quantum groups) as the underlying geometry when equipped with left and bi-covariant differential calculi in the sense of \cite{Wor:dif}. One innovation is that we will consider also  calculi that are left invariant but only partially right invariant.

Section 5 then applies the theory to finite groups $G$ with algebra $C(G)$ of functions and bicovariant calculus of differential forms, as considered particularly in \cite{BMHDSfinitegrps} and \cite{Ma:rief}.  We focus calculations on $S_3$ with its standard 3D bicovariant calculus, giving a full analysis of the moduli space of connections under the various conditions of the general theory (but we limit ourselves to the standard Euclidean metric). Briefly, we find a 1-parameter moduli of  torsion free cotorsion free $*$-compatible `generalised Levi-Civita connections'; and we  find a discrete (finite) moduli of metric preserving and $*$-compatible connections, but with torsion. Thus $*$-compatibility requires us either to only weakly preserve the metric or to have torsion as a necessary feature of the discrete Riemannian geometry. This contrasts with the unique torsion free and cotorsion free connection previously found in \cite{Ma:rief}, which turns out to be not $*$-compatible. Precisely four  of our metric perserving $*$-compatible connections have $\sigma$ a braiding. One of these is the left Maurer-Cartan connection which we characterise now in these geometric terms, while the others are   `nonstandard' variants of it. The braiding obeys $\sigma^3=\pm \id$,  see Corollary~\ref{corbras3}.

Finally, the critical case of the quantum group $C_q[SU_2]$ is fully analysed in Section~6 and applied to the quantum sphere $C_q[S^2]$ as an example of the theory with nontrivial cotangent bundle. On the quantum group we use the standard  Woronowicz 3D calculus and for practical reasons we limit ourselves to a three parameter family of diagonal metrics $g=(g_{++},g_{00},g_{--})$ characterised as left invariant and right $U(1)$-invariant, and to similarly covariant connections $\nabla$. This contrasts with previous work \cite{Ma:rief} where the 4D calculus and only the Killing metric were considered. The theory turns out to work well and provide a unique `Levi-Civita'  connection for each point of the moduli space of metrics.

This result, see Theorem~\ref{nceicnpncuu}, pertains to the existence of a unique torsion free, metric and $*$-compatible connection with classical limit as $q\to 1$.  If we drop the last requirement we have three other `purely quantum' connections fulfilling the remaining requirements. This ambiguity is similar to that encountered in the classification of the differential calculi themselves where, typically, one has `purely quantum variants' of the  canonical choice obtained by twisting by discrete characters. We also look at when $\sigma$ obey the braid relations and find a unique such connection that is metric and $*$-compatible, but necessarily with torsion, Corollary~\ref{kzdvvcdijsbc}.
It is this connection or its `purely quantum'  counterpart, which turns out to project down to a connection on the quantum sphere where we find that it coincides with the cannonical $q$-deformed Levi-Civita connection found in \cite{Ma:rsph} as we see in Section~6.2.

In both cases, our calculations make extensive use of computer algebra (Mathematica in our case). Of course, as anyone who has worked on General Relativity will know, finding explicit Riemannian metrics is not necessarily easy in the classical case, but we expect that more general results will be needed to facilitate handling more complicated noncommutative geometries.

The paper was mainly written during the visit July-December 2006 of the authors to the  Isaac Newton Institute. We thank the institute for their support.

\section{Covariant derivatives, duals and the bar functor}

\subsection{Bimodule covariant derivatives}

Let $A$ be a unital algebra over $\C$. The basic algebraic notions of course work more generally. Suppose that the algebra $A$ has a differential structure $(\Omega A,\extd)$ in the sense of a differential graded exterior algebra $\Omega A=\oplus_n\Omega^n A$ with $d$ increasing degree and obeying a graded Leibniz rule and $\extd^2=0$. We suppose that $\Omega^1$ generates the exterior algebra over $A$. The notion of a covariant derivative in this context is standard:

\begin{defin} Given a left $A$-module $E$, a left $A$-covariant derivative
is a map $\nabla:E\to\Omega^{1}A\tens_{A}E$ which obeys the condition
$\nabla(a.e)=da\tens e+a.\nabla e$
for all $e\in E$ and $a\in A$.
\end{defin}

In classical differential geometry there is no difference in whether we multiply a section by a function on the left or right. In the noncommutative case there is a difference, and
for a bimodule we could require the Leibnitz rule for both
left and right multiplication, but this would turn out to be too restrictive. Instead,
following \cite{DVMic,DVMass} and \cite{Mourad}, we make the following definition:

\begin{defin} \label{ppll}  A bimodule covariant derivative on an
 $A$-bimodule $E$ is a triple $(E,\nabla,\sigma)$,
where $\nabla:E\to\Omega^{1}A\tens_{A}E$
is a left $A$-covariant derivative, and $\sigma:E\tens_A\Omega^1 A\to
\Omega^1 A \tens_A E$ is a bimodule map obeying
\[
\nabla(e.a)\,=\,\nabla(e).a\,+\,\sigma(e\tens da)\ ,\qquad \forall\, e\in E,\ a\in A
\]
\end{defin}

Now we consider one of the most immediate reasons to define the
bimodule covariant derivative, that is to have a covariant derivative on tensor
products of bimodules. As mentioned in the introduction, this occurred in 
\cite{BMHDSfinitegrps}.

\begin{propos} \label{tens45}  Given $(E,\nabla_E,\sigma_E)$ a
bimodule covariant derivative on the bimodule $E$ and $\nabla_F$ a left
covariant derivative on the left module $F$, there is a left
$A$-covariant derivative on $E\tens_{A}F$ given by
\begin{eqnarray*}
\nabla_{E\tens F} &=& \nabla_E\tens
\id_{F}+(\sigma_E\tens\id_{F})(\id_{E}\tens\nabla_F)
\end{eqnarray*}
Further if $F$ is also an $A$-bimodule with a bimodule covariant
derivative $(\nabla_F,\sigma_F)$, then there is a
bimodule covariant
derivative $(\nabla_{E\tens_A F},\sigma_{E\tens_A F})$ on $E\tens_A F$ with
\begin{eqnarray*}
\sigma_{E\tens F} \,=\, (\sigma_{E}\tens
\id)(\id\tens \sigma_{F})\ .
\end{eqnarray*}
\end{propos}
\proof Applying $\nabla_{E\tens F}$ to $e\tens a.f$ (for $a\in A$) we get
\begin{eqnarray*}
\nabla_E e\tens a.f+(\sigma_E\tens\id_{F})(e\tens \extd a\tens f+e\tens a.\nabla f)
\end{eqnarray*}
Aplying the formula to $e.a\tens f$ we get $\nabla_E(e.a)\tens f+
(\sigma_E\tens \id_{F})(e.a\tens \nabla_F f)$, and these are the same by
definition of $\sigma_E$.
 This shows that the given formula is
well defined on $E\tens_{A}F$.  The left multiplication property is
true because $\sigma_E$ is a left $A$-module map.

For the second part,
\begin{eqnarray*}
    \sigma_{E\tens F}(e\tens f\tens a.db) &=& \nabla(e\tens f.ab)-
    \nabla(e\tens f.a).b \cr
    &=& \nabla(e)\tens f.ab + (\sigma_{E}\tens\id )(e\tens \nabla(f.ab)) \cr
    &&- \nabla(e)\tens f.ab-(\sigma_{E}\tens\id )(e\tens
    \nabla(f.a)b)\ .\quad\square
\end{eqnarray*}

\begin{defin} \label{skudg}
The category ${}_A\mathcal{E}_A$ consists of objects $A$-bimodule covariant derivatives
 $(E,\nabla,\sigma)$ where $\sigma:E\tens_A \Omega^1 A \to  \Omega^1 A \tens_A E$
is invertible.
 The morphisms are bimodule maps $\theta:E\to F$ which are preserved by the
covariant derivatives, i.e.\
\begin{eqnarray*}
\nabla\circ\theta=(\id\tens\theta)\nabla:E\to\Omega^1 A\tens_A F\ .
\end{eqnarray*}
Then \ref{tens45} makes  ${}_A\mathcal{E}_A$ into a monoidal category.
The identity for the tensor product is the bimodule $A$, with
$\nabla_A=\extd:A\to \Omega^1 A\tens_A A=\Omega^1 A$, and
$\sigma_A$ the is identity map $A\tens_A \Omega^1 A$ to
$\Omega^1 A \tens_A A$ when both sides are identified with $\Omega^1A$.
\end{defin}

A monoidal category here means a category with a $\tens$ functor obeying some standard axioms of associativity. In this paper all our monoidal categories will  be `built' on the category ${}_A\CM_A$ of bimodules over $A$ equipped with further structure such as above. Here $\tens$ on the underlying bimodules is $\tens_A$ and morphisms are among other things bimodule maps (in other words, there is a forgetful functor to ${}_A\CM_A$). It may surprise the reader that we do not impose in ${}_A\CE_A$ any conditions between a morphism $\theta$ and $\sigma$.
In fact there is an equation which is automatically true:

\begin{propos} \label{bchdjsabv}
Suppose that $\theta: (E,\nabla_E,\sigma_E)\to (F,\nabla_F,\sigma_F)$ is a morphism
in ${}_A\mathcal{E}_A$ as described in \ref{skudg}. Then
\begin{eqnarray*}
(\id\tens\theta)\, \sigma_E\,=\, \sigma_F\, (\theta\tens\id) : E\tens_A\Omega^1 A
\to \Omega^1 A \tens_A F\ .
\end{eqnarray*}
\end{propos}
\proof As $\sigma_E$, $\sigma_F$ and $\theta$ are all
bimodule maps, it is enough to check the equation on $e\tens\extd a$:
\begin{eqnarray*}
(\id\tens\theta)\, \sigma_E(e\tens\extd a) &=& (\id\tens\theta)
(\nabla_E(e.a)-\nabla_E(e).a) \cr
&=& \nabla_F(\theta(e).a)-\nabla_F(\theta(e)).a \cr
&=& \sigma_F(\theta(e)\tens\extd a)\ .\quad\square
\end{eqnarray*}

\medskip Note that the condition in \ref{bchdjsabv} can be viewed as a rather weaker notion of morphism between objects of ${}_A\mathcal{E}_A$, again giving a bar category.

\subsection{Torsion}
Here we shall restrict attention to a left covariant derivative $\nabla$ on $\Omega^1 A$.
Following \cite{Ma:rief}  make the following definition, which
 in the classical case reduces to the usual one:

\begin{defin}
The torsion of a left $A$-covariant derivative $\nabla$ on $\Omega^1 A$ is the left $A$-module map $T=\wedge\nabla-\extd:\Omega^1 A \to \Omega^2 A$.
\end{defin}

That it is a left module map follows easily from the definition of a covariant derivative
\begin{eqnarray*}
T(a.\xi) &=& \wedge(a.\nabla\xi)+\wedge(\extd a\tens\xi) -a.\extd\xi
+\extd a\wedge\xi= a.T(\xi)\
\end{eqnarray*}
for all $\xi\in\Omega^1 A$ and $a\in A$. We
now apply the definition of torsion to a bimodule covariant derivative, and obtain
the following result:

\begin{propos} \label{khshdvbcskhjv}
 The torsion of a
bimodule covariant derivative $(\Omega^1 A,\nabla,\sigma)$
 is a bimodule map if and only if
 \begin{eqnarray*}
\mathrm{image}(\id+\sigma)\subset \ker(\wedge:\Omega^1 A \tens_A
\Omega^1 A\to \Omega^2 A)\ .
\end{eqnarray*}
We say in this case that $\nabla$ is {\em torsion-compatible}.
\end{propos}
\proof For  $\xi\in\Omega^1 A$ and $a\in A$,
\begin{eqnarray*}
T(\xi.a) &=& \wedge\nabla\xi.a+\wedge\sigma(\xi\tens\extd a) -\extd\xi.a
+\xi\wedge\extd a \cr
&=& T(\xi).a+\wedge\sigma(\xi\tens\extd a)
+\xi\wedge\extd a\ .
\end{eqnarray*}
Now use the fact that both $\wedge$ and $\sigma$ are bimodule maps,
and the density condition.
\quad$\square$

\medskip Note that the most important case classically is where the torsion vanishes, and
as zero is a bimodule map, in this case \ref{khshdvbcskhjv} applies.

\subsection{Finitely generated projective modules and duals}
An object $X$ in a tensor category has a left dual if there is an object
$X'$ and morphisms $\ev^L_X:X'\tens X\to 1_{\mathcal{C}}$
(evaluation) and
$\coev^L_X:1_{\mathcal{C}}\to X\tens X'$ (coevaluation) so that
\begin{eqnarray*}
l^{-1}_X(\id\tens\ev_X)\Phi(\coev_X\tens \id)r_X  &=& \id_X:X\to X\ , \cr
r_{X'}^{-1}(\ev_X\tens\id)\Phi^{-1}(\id\tens\coev_{X})l_{X'}  &=& \id_{X'}:X'\to X'\ .
\end{eqnarray*}
Here we use the canonical isomorphisms $r_X:X\to 1_{\mathcal{C}} \tens X$ and
$l_X:X\to X\tens 1_{\mathcal{C}} $, which in our bimodule category case are trivial.
It is natural to ask when an object $(E,\nabla_E,\sigma_E)$ in ${}_A\mathcal{E}_A$ has a left dual,
but first we forget about covariant derivatives, and just consider bimodules.

Suppose that the bimodule $E$ is finitely generated projective as a right $A$-module.
We shall define the dual $E'=\mathrm{Hom}_A(E,A)$, the right
module maps from $E$ to $A$. Remember that a
 right $A$-module $E$ is said to be finitely
    generated projective if
there are $e^i\in E$ and $e_i\in E'$ (for integer $1\le i\le n$)
(the `dual basis') so that for all $f\in E$, $f=\sum e^i.e_i(f)$.
From this it follows directly that $\alpha=\sum \alpha(e^i).e_i$
for all $\alpha\in E'$.
We will use implicit summation over the dual basis index. Now the evaluation map
$\ev^L_E:E'\tens E\to A$ is just $\ev^L_E(\alpha\tens e)=\alpha(e)$, and the coevaluation
$\coev^L_E:A\to E\tens E'$ is given by $\coev^L_E(1)=\sum e^i\tens e_i$.
To show uniqueness at various points in what follows it will be useful to use the following result:

\begin{lemma}\label{endo} If an $A$-bimodule $E$ is right finitely generated projective,
and $F$ is a right $M$-module,
 there is an isomorphism $\vartheta:F\tens_A E' \to {\rm Hom}_A(E,F)$
defined by $\vartheta(f\tens\alpha)(e)=f.\alpha(e)$.
\end{lemma}
\proof The inverse map is $\vartheta^{-1}(T)=\sum T(e^i)\tens e_i$.\quad$\square$

\begin{propos} \label{camab} Given $(E,\nabla_E,\sigma_E)$ in ${}_A\mathcal{E}_A$,
where $E$ is finitely generated projective as a right $A$-module,
there is a unique bimodule covariant derivative
$(E',\nabla_{E'},\sigma_{E'})$ on $E'$ so that the map $\ev_E^L:E'\tens_A E\to
A$ is a morphism in ${}_A\CE_A$.
(Here we take the identity $(A,\extd,\id)$ in ${}_A\CE_A$.) It is
defined in terms of the dual basis of $E$ by
\begin{eqnarray*}
\sigma_{E'}(\alpha\tens\xi) &=&
 (\ev\tens\id)(\id\tens\sigma_E^{-1})(\alpha\tens\xi\tens e^i) \tens e_i\ ,
\cr
\nabla_{E'}\alpha &=&  \extd(\alpha(e^i))\tens e_i\,-\,
(\ev\tens\id)(\id\tens\sigma_E^{-1}\nabla_E)(\alpha\tens e^i)
\tens e_i\ .
\end{eqnarray*}
\end{propos}
\proof First we check that the formulae give a left covariant derivative.
For $a\in A$ and $\alpha\in E'$,
\begin{eqnarray*}
\nabla_{E'}(a.\alpha) &=&  \extd(a.\alpha(e^i))\tens e_i-(\ev\tens\id)
(\id\tens\sigma_E^{-1})(\id\tens\nabla_E)(a.\alpha\tens e^i)
\tens e_i \cr
&=& \extd a\tens \alpha(e^i).e_i \,+\, a.\nabla_{E'}(\alpha)\,=\,
\extd a\tens \alpha \,+\, a.\nabla_{E'}(\alpha)\ .
\end{eqnarray*}
To see that the covariant derivatives preserve the evaluation:
\begin{eqnarray}\label{nabla-eval}
(\id\tens\ev)(\nabla_{E'}(\alpha)\tens e) &=&
\extd(\alpha(e_i)).\alpha_i(e) \cr && -\,
(\ev\tens\id)(\id\tens\sigma_E^{-1})(\id\tens\nabla_E)
(\alpha\tens e_i).\alpha_i(e) \cr
&=&  \extd(\alpha(e))
\,-\,(\ev\tens\id)(\id\tens\sigma_E^{-1})(\id\tens\nabla_E)
(\alpha\tens e) \cr
&& -\, \alpha(e_i)\,\extd(\alpha_i(e))+
(\ev\tens\id)(\alpha\tens e_i\tens \extd(\alpha_i(e)))\cr
&=&
\extd(\alpha(e))
-(\id\tens\ev)(\sigma_{E'}\tens\id)(\id\tens\nabla_E)(\alpha\tens e)\ .
\end{eqnarray}
The $\nabla_{E'}$ with this property (\ref{nabla-eval}) is unique by \ref{endo}.
Finally we check the compatibility condition in \ref{ppll}, using (\ref{nabla-eval}):
\begin{eqnarray*}
(\id\tens\ev)(\nabla_{E'}(\alpha.a)\tens e) &=&
  \extd((\alpha.a)(e))
-(\ev\tens\id)(\id\tens\sigma_E^{-1}\nabla_E)(\alpha.a\tens e) \cr
 &=&
  \extd(\alpha(a.e))
-(\ev\tens\id)(\id\tens\sigma_E^{-1}\nabla_E)(\alpha\tens a.e) \cr
&&+\,(\ev\tens\id)(\id\tens\sigma_E^{-1})(\alpha\tens \extd a\tens e) \cr
&=& (\id\tens\ev)(\nabla_{E'}(\alpha).a\tens e)
\cr
&& +\, (\id\tens\ev)(\sigma_{E'}\tens\id)(\alpha\tens \extd a\tens e)\ .
\quad\square
\end{eqnarray*}

\begin{lemma}\label{kronprop}
 Given $(E,\nabla_E,\sigma_E)$ in ${}_A\mathcal{E}_A$,
where $E$ is finitely generated projective as a right $A$-module,
the Kroneker delta, $\delta_E=\coev_E^L(1)
=\sum e^i\tens e_i$
has the following properties:

a)\quad $a.\delta_E=\delta_E.a$ for all $a\in A$.

b)\quad $\nabla_{E\tens E'}(\delta_E)=0$.

c)\quad $\sigma_{E\tens E'}
(\delta_E\tens\xi)=\xi\tens\delta_E$ for all $\xi\in\Omega^1 M$.

\end{lemma}
\proof  By \ref{endo}, to prove (a) we only have to show that
\begin{eqnarray}\label{test1}
a.e &=& (\id\tens\ev)(\delta_E.a\tens e)\ .
\end{eqnarray}
But the right hand side of (\ref{test1}) is
\[
 e^i.(e_i.a)(e)\,=\, e^i.e_i(a.e)\,=\,a.e
\]
for all $e\in E$, as required.

By \ref{endo}, to prove (b) we only have to show that
$(\id^2\tens\ev)(\nabla_{E\tens E'}\delta_E\tens e)=0$ for all $e\in E$. Then, using (\ref{nabla-eval}),
\begin{eqnarray*}
(\id^2\tens\ev)(\nabla\delta_E\tens e) &=& (\id^2\tens\ev)\Big(\nabla e^i\tens
e_i\tens e
+(\sigma_E\tens\id^2)(e^i\tens\nabla e_i\tens e)\Big) \cr
&=& (\nabla e^i).e_i(e)\,+\,\sigma_E(e^i\tens(\id\tens\ev)(\nabla e_i\tens e))\cr
&=& (\nabla e^i).e_i(e) \,+\, \sigma_E(e^i\tens \extd e_i(e))\cr &&-\,
\sigma_E(e^i\tens(\id\tens\ev)(\sigma_{E^*}\tens\id)(e_i\tens\nabla e)) \cr
&=&\nabla(e^i.e_i(e))-
\sigma_E(e^i\tens(\id\tens\ev)(\sigma_{E^*}\tens\id)(e_i\tens\nabla e)) \cr
&=& \nabla(e)
- \sigma_E(e^i\tens(\ev\tens\id)(\id\tens\sigma_{E}^{-1})(e_i\tens\nabla e))\ .
\end{eqnarray*}
Now substitute $\sigma_{E}^{-1}\nabla e= f_j\tens\eta_j\in E\tens\Omega^1 M$
(summation implicit), giving
\begin{eqnarray*}
(\id^2\tens\ev)(\nabla\delta_E\tens e)&=& \nabla(e)\,-\,
\sigma_E(e^i\tens e_i(f_j).\eta_j  ) \cr
&=& \nabla(e)\,-\,
\sigma_E(e^i.e_i(f_j)\tens\eta_j  ) \cr
&=& \nabla(e)\,-\,
\sigma_E(f_j\tens\eta_j  ) \,=\,0
\end{eqnarray*}
Now applying $\nabla_{E\tens E'}$ to $\delta_E.a=a.\delta_E$ and using (b) gives (c).
\quad$\square$

We shall also find it convenient to consider the right dual $X^\circ={}_A\Hom(X,A)$ in the case where $X$ is a finitely generated projective left $A$-module. The corresponding
evaluation and coevaluation maps will be written
 $\ev^R_X:X\tens X^\circ\to 1_{\mathcal{C}}$
 and
$\coev^R_X:1_{\mathcal{C}}\to X^\circ\tens X$.
Now we have a \textit{right} connection on $E^\circ$, where we write $\coev^R_X(1)=f_i\tens f^i$,
\begin{eqnarray*}
\hat\nabla_{X^\circ}\alpha &=& f_i\tens (\extd(\ev^R(f^i\tens\alpha))-(\id \tens\ev^R)(\nabla_X(f^i)\tens\alpha))
\end{eqnarray*}
The reader can now check that
\begin{eqnarray*}
\extd(\ev^R(e\tens\alpha)) &=& (\id\tens\ev^R)(\nabla e\tens\alpha) + (\ev^R\tens\id)(e\tens\hat\nabla\alpha)\ .
\end{eqnarray*}

\begin{eqnarray*}
(\ev^R\tens\id)(e\tens\hat\nabla(a.\alpha)) &=&
\extd(\ev^R(e\tens a.\alpha)) - (\id\tens\ev^R)(\nabla e\tens a.\alpha) \cr
 &=&
\extd(\ev^R(e.a\tens \alpha)) - (\id\tens\ev^R)((\nabla e).a\tens \alpha) \cr
 &=&
\extd(\ev^R(e.a\tens \alpha)) - (\id\tens\ev^R)((\nabla e.a)\tens \alpha) \cr
&& +\ (\id\tens\ev^R)(\sigma(e\tens\extd a)\tens\alpha) \cr
 &=&
(\ev^R\tens\id)(e.a\tens\hat\nabla(\alpha))  + (\id\tens\ev^R)(\sigma(e\tens\extd a)\tens\alpha)\ ,
\end{eqnarray*}
so
\begin{eqnarray*}
(\ev^R\tens\id)(e\tens(\hat\nabla(a.\alpha)-a.\hat\nabla(\alpha))) &=&
(\id\tens\ev^R)(\sigma(e\tens\extd a)\tens\alpha)\ .
\end{eqnarray*}
Then
\begin{eqnarray*}
\hat\nabla(a.\alpha)-a.\hat\nabla(\alpha) &=&
f_i\,(\ev^R\tens\id)(f^i\tens(\hat\nabla(a.\alpha)-a.\hat\nabla(\alpha))) \cr
 &=&
f_i\,(\id\tens\ev^R)(\sigma_E(f^i\tens\extd a)\tens\alpha)\ .
\end{eqnarray*}
We will label
\begin{eqnarray*}
 \sigma^{-1}_{E^\circ} (\xi\tens \alpha) &=&
f_i\,(\id\tens\ev^R)(\sigma_E(f^i\tens\xi)\tens\alpha)\ ,
\end{eqnarray*}
and assume that it is invertible with inverse $ \sigma_{E^\circ}$. Then $\nabla_{E^\circ}=
\sigma_{E^\circ}\, \hat\nabla$ is a left connection on $E^\circ$ with
\begin{eqnarray*}
\nabla_{E^\circ}(\alpha.a) \,=\, \nabla_{E^\circ}(\alpha).a+ \sigma_{E^\circ}(\alpha\tens\extd a)\ .
\end{eqnarray*}

\subsection{Covariant derivatives and the bar functor} \label{ikcvsaucfv}
Now we assume that $A$ is a star algebra. Then the category ${}_A\mathcal{M}_A$
of $A$-bimodules becomes a bar category (see \cite{barcats}). This is a monoidal
category equipped with a functor bar that sends every object $E$ to a `conjugate' one
$\overline{E}$ in a manner that reverses order up to a natural equivalence $\Upsilon$. There are some auxiliary data and constraints such as natural isomorphisms ${\mathrm{bb}}_E:E\to \overline{\overline{E}}$ needed to make the theory work. In ${}_A\CM_A$ and the morphisms $\Upsilon$ and the other data will all be the `obvious' maps given by flip and the application of $*$ etc, so the bar category structure is not essential but remains a way of thinking clearly about the constructions. If $E$ is a bimodule then $\overline{E}$ is identified with $E$ as a set but has the conjugate actions $a.\bar e=\overline{e.a^*}$ and $\bar e.a=\overline{a^*.e}$. Here $\bar e$ denotes $e\in E$ viewed in $\overline{E}$.
In a bar category we define a star object to be an object equipped with a
morphism $\star:E\to \overline{E}$ so that $\bar\star\star(e)=\overline{\overline{e}}$,
where we use the notation $\star(e)=\overline{e^*}$.

We suppose that $\Omega^1 A$ is a star object in ${}_A\CM_A$ in such a manner that
$\star\,\extd=\overline{\extd}\,\star:A\to\overline{\Omega^1 A}$.
This comes down to the concrete assumptions that $\Omega^1A$ has an antilinear involution $\star$ compatible with $\extd$ in the sense
 \[ \star(a.\extd b)=\overline{\extd b^*.a^*}\]
 If we define $(\extd a)^*=\extd a^*$ as usual in noncommutative geometry, then the above condition amounts to $(a.\extd b)^*=\extd b^*.a^*$, i.e.\ the lowest degree part of the usual notion of a star differential graded algebra.

 Now we come to the main reason why we have to distinguish between
 $E$ and
$\overline{E}$ in an obvious manner, because the covariant derivatives on them look very different.

\begin{propos}
Given a left bimodule connection $(\nabla_E,\sigma_E)$ on $E$,
there is a right connection $\hat\nabla$ on $\overline{E}$ given by
\begin{eqnarray*}
\hat\nabla\, \overline e\,=\, (\id\tens\star^{-1})\Upsilon(\overline{\nabla_E e})\ .
\end{eqnarray*}
Further the map $\hat\sigma=(\id\tens\star^{-1})\Upsilon\,\overline{\sigma_E}\,
\Upsilon^{-1}(\star\tens\id)$ satisfies
\begin{eqnarray*}
\hat\nabla(a.\overline{e}) \,=\, a.\hat\nabla(\overline{e}) + \hat\sigma(\extd a\tens \overline{e}) \ .
\end{eqnarray*}
\end{propos}
\proof First we check the right Liebnitz property:
\begin{eqnarray*}
\hat\nabla(\overline{e}.a) &=& \hat\nabla(\overline{a^*.e}) \cr
&=& (\id\tens\star^{-1})\Upsilon(\overline{\nabla_E (a^*.e)}) \cr
&=& (\id\tens\star^{-1})\Upsilon(\overline{\extd a^*\tens e +a^*.\nabla_E (e)}) \cr
&=& \overline{e}\tens\extd a + (\id\tens\star^{-1})\Upsilon(\overline{a^*.\nabla_E (e)}) \cr
&=& \overline{e}\tens\extd a + (\id\tens\star^{-1})\Upsilon(\overline{\nabla_E (e)}).a\ .
\end{eqnarray*}
Now we check the property for $\hat\sigma$:
\begin{eqnarray*}
\hat\nabla(a.\overline{e}) &=& \hat\nabla(\overline{e.a^*}) \cr
&=&  (\id\tens\star^{-1})\Upsilon(\overline{\nabla_E (e.a^*)}) \cr
&=&  (\id\tens\star^{-1})\Upsilon(\overline{\nabla_E (e).a^*})+
 (\id\tens\star^{-1})\Upsilon\,\overline{\sigma_E(e\tens\extd a^*)}\ .\quad\square
\end{eqnarray*}

\medskip
On the assumption that $\sigma_E$ is invertible,
we can apply $\hat\sigma^{-1}$ to $\hat\nabla$ to get the following:

\begin{theorem}\label{sliogcb}
Suppose that $A$ is a
 star algebra which has a differential structure $(\Omega^1 A,\extd)$
 so that $\Omega^1 A$ is a star object and
 $\star\,\extd=\overline{\extd}\,\star:A\to\overline{\Omega^1 A}$.
 Then  ${}_A\mathcal{E}_A$
described in \ref{skudg} is a bar category with
$\overline{(E,\nabla_E,\sigma_E)}=(\overline{E},\nabla_{\bar E},\sigma_{\bar E})$ given by
\begin{eqnarray*}
\nabla_{\overline{E}}(\overline{e}) &=& (\star^{-1}\tens\id) \Upsilon\,\overline{\sigma_E^{-1}
\nabla_E(e)}\ ,\cr
\sigma_{\overline{E}} &=&  (\star^{-1}\tens\id) \Upsilon\, \overline{\sigma_E^{-1}}\, \Upsilon^{-1}
(\id\tens\star)\ .
\end{eqnarray*}
\end{theorem}
\proof This is mostly inherited from the bar category structure of ${}_A\mathcal{M}_A$. What remains to check is that
various things are objects or morphisms in ${}_A\mathcal{E}_A$, and
this is contained in the lemmas \ref{fghgclem2},
\ref{fghgclem3} and \ref{fghgclem4}.\quad$\square$

\begin{lemma}  \label{fghgclem2}
Given a morphism $\theta:(E,\nabla_E,\sigma_E)\to
(F,\nabla_F,\sigma_F)$ in ${}_A\mathcal{E}_A$,
then the usual formula for $\overline{\theta}$
(i.e.\ $\overline{\theta}(\overline{e})=\overline{\theta(e)}$)
gives a morphism
$\overline{\theta}:(\overline{E},\nabla_{\bar E},\sigma_{\bar E}) \to
(\overline{F},\nabla_{\bar F},\sigma_{\bar F})$.
\end{lemma}
\proof We just have to check that the bimodule map
$\overline{\theta}:\overline{E}\to \overline{F}$ is preserved by the corresponding
covariant derivatives.
\begin{eqnarray*}
\nabla_{\bar F}(\overline{\theta}(\overline{e})) &=&
(\star^{-1}\tens \id)\Upsilon_{F,\Omega^1 A}
(\overline{\sigma_F^{-1}\nabla_E(\theta(e))}) \cr
 &=&
(\star^{-1}\tens \id)\Upsilon_{F,\Omega^1 A}
(\overline{\sigma_F^{-1}(\id\tens\theta)\nabla_E(e)}) \cr
 &=&
(\star^{-1}\tens \id)\Upsilon_{F,\Omega^1 A}
(\overline{(\theta\tens\id)\sigma_E^{-1}\nabla_E(e)})\cr
 &=&
(\star^{-1}\tens \overline{\theta})\Upsilon_{E,\Omega^1 A}
(\overline{\sigma_E^{-1}\nabla_E(e)}) \cr
&=& (\id\tens \overline{\theta})\,\nabla_{\bar E}(\overline{e})\ .\quad\square
\end{eqnarray*}

\begin{lemma} \label{fghgclem3} Given the operations defined in \ref{sliogcb},
for objects $(E,\nabla_E,\sigma_E)$ and
$(F,\nabla_F,\sigma_F)$ in ${}_A\mathcal{E}_A$, the bimodule map
$\Upsilon_{F,E}:\overline{F\tens E}\to \overline{E}\tens_A \overline{F}$
 (defined just as for
${}_A\mathcal{M}_A$) is a morphism in ${}_A\mathcal{E}_A$.
\end{lemma}
\proof We need to verify the following equation:
\begin{eqnarray}\label{skjdhjcvbver}
(\id\tens\Upsilon_{F,E})
\nabla_{\overline{F\tens E}} \,=\, \nabla_{\bar E \tens\bar F}\,\Upsilon_{F,E}\ .
\end{eqnarray}
Begin with
\begin{eqnarray*}
(\id\tens\Upsilon_{F,E})
\nabla_{\overline{F\tens E}}(\overline{f\tens_A e}) &=&
(\star^{-1}\tens \Upsilon_{F,E})\Upsilon_{F\tens E,\Omega^1 A}
(\overline{\sigma_{F\tens E}^{-1}\nabla_{F\tens E}(f\tens_A e)})\ .
\end{eqnarray*}
On the other hand,
\begin{eqnarray*}
\nabla_{\bar E \tens\bar F}\Upsilon_{F,E}(\overline{f\tens_A e}) &=&
\nabla_{\bar E \tens\bar F}(\bar e\tens_A\bar f) \cr
 &=& \nabla_{\bar E}\bar e\tens_A\bar f
+(\sigma_{\bar E}\tens\id)(\bar e\tens_A\nabla_{\bar F}\bar f) \cr
&=& (\star^{-1}\tens \id)\Upsilon_{E,\Omega^1 A}(\overline{\sigma_E^{-1}\nabla_E(e)}) \tens_A
\bar f \cr
&&+\, \big((\star^{-1}\tens \id)\Upsilon_{E,\Omega^1 A}
\overline{\sigma_E^{-1}}\Upsilon_{\Omega^1 A,E}^{-1}(\id\tens *)\tens\id\big)
(\bar e\tens_A\nabla_{\bar F}\bar f) \ .
\end{eqnarray*}
Then the statement we are asked to verify becomes
\begin{eqnarray} \label{hdvfszdhv}
&&(\id\tens \Upsilon_{F,E})\Upsilon_{F\tens E,\Omega^1 A}
(\overline{\sigma_{F\tens E}^{-1}\nabla_{F\tens E}(f\tens_A e)})\cr
&=& (\Upsilon_{E,\Omega^1 A}\, \overline{\sigma_E^{-1}}\tens\id)
[\overline{\nabla_E(e)} \tens_A\bar f +\big(\Upsilon_{\Omega^1 A,E}^{-1}\tens\id\big)
(\bar e\tens_A(*\tens\id)\nabla_{\bar F}\bar f) ]\ .
\end{eqnarray}
Now we have, by direct calculation,
\begin{eqnarray*}
(\Upsilon_{E,\Omega^1 A}^{-1}\tens\id)(\id\tens \Upsilon_{F,E})\Upsilon_{F\tens E,\Omega^1 A}
\,=\,\Upsilon_{F,E\tens\Omega^1 A}\ ,
\end{eqnarray*}
and substituting this into
 (\ref{hdvfszdhv}) leaves us with having to verify
 \begin{eqnarray} \label{hdvfszdhv888}
&&\Upsilon_{F,E\tens\Omega^1 A}
(\overline{\sigma_{F\tens E}^{-1}\nabla_{F\tens E}(f\tens_A e)})\cr
&=& ( \overline{\sigma_E^{-1}}\tens\id)
[\overline{\nabla_E(e)} \tens_A\bar f +\big(\Upsilon_{\Omega^1 A,E}^{-1}\tens\id\big)
(\bar e\tens_A(*\tens\id)\nabla_{\bar F}\bar f) ]\ .
\end{eqnarray}
Here the expression under the large bar is
\begin{eqnarray*}
\sigma_{F\tens E}^{-1}\nabla_{F\tens E}(f\tens_A e) &=&
(f\tens_A\sigma_E^{-1} \nabla_E e)+
(\id\tens\sigma_E^{-1})(\sigma_F^{-1} \nabla_F f\tens_A e)\ .
\end{eqnarray*}
As
 \begin{eqnarray*}
\Upsilon_{F,E\tens\Omega^1 A}
(\overline{f\tens_A\sigma_E^{-1} \nabla_E e})
&=& ( \overline{\sigma_E^{-1}}\tens\id)
[\overline{\nabla_E(e)} \tens_A\bar f ]\ ,
\end{eqnarray*}
we can again simplify the required condition to
 \begin{eqnarray*}
&&\Upsilon_{F,E\tens\Omega^1 A}
(\overline{(\id\tens\sigma_E^{-1})(\sigma_F^{-1} \nabla_F f\tens_A e)})\cr
&=& ( \overline{\sigma_E^{-1}}\tens\id)
[\big(\Upsilon_{\Omega^1 A,E}^{-1}\tens\id\big)
(\bar e\tens_A(*\tens\id)\nabla_{\bar F}\bar f) ]\ ,
\end{eqnarray*}
 which immediately undergoes cancellation to become
  \begin{eqnarray*}
\Upsilon_{F,\Omega^1 A\tens E}
(\overline{\sigma_F^{-1} \nabla_F f\tens_A e})
&=&
\big(\Upsilon_{\Omega^1 A,E}^{-1}\tens\id\big)
(\bar e\tens_A(*\tens\id)\nabla_{\bar F}\bar f) \ .
\end{eqnarray*}
 On using the definition of $\nabla_{\bar F}$ this condition becomes
   \begin{eqnarray} \label{hdvfszdhv88899}
\Upsilon_{F,\Omega^1 A\tens E}
(\overline{\sigma_F^{-1} \nabla_F f\tens_A e})
&=&
\big(\Upsilon_{\Omega^1 A,E}^{-1}\tens\id\big)
(\bar e\tens_A \Upsilon_{F,\Omega^1 A} \overline{\sigma_F^{-1}\nabla_{F} f}) \ .
\end{eqnarray}
By direct calculation, we have
    \begin{eqnarray*}
(\id\tens \Upsilon_{F,\Omega^1 A}^{-1})
(\Upsilon_{\Omega^1 A,E}\tens\id)\Upsilon_{F,\Omega^1 A\tens E}
&=&
 \Upsilon_{F\tens\Omega^1 A,E}  \ ,
\end{eqnarray*}
 and using this in (\ref{hdvfszdhv88899}) gives
    \begin{eqnarray*}
     \Upsilon_{F\tens\Omega^1 A,E}
(\overline{\sigma_F^{-1} \nabla_F f\tens_A e})
&=&
\bar e\tens_A  \overline{\sigma_F^{-1}\nabla_{F} f} \ .\quad\square
\end{eqnarray*}

\begin{lemma} \label{fghgclem4}
 Given the operations defined in \ref{sliogcb},
for an object $(E,\nabla_E,\sigma_E)$ in ${}_A\mathcal{E}_A$, the bimodule map
$\mathrm{bb}_E:E\to \overline{\overline{E}} $
 (defined just as for
${}_A\mathcal{M}_A$) is a morphism in ${}_A\mathcal{E}_A$.
\end{lemma}
\proof We need to verify the following equation:
\begin{eqnarray*}
\nabla_{\overline{\bar E}} \circ\mathrm{bb}_E \,=\,(\id\tens \mathrm{bb}_E)\nabla_E\ .
\end{eqnarray*}
Begin with
 \begin{eqnarray*}
\nabla_{\overline{\bar E}}(\overline{\bar e})
&=& (\star^{-1}\tens \id)\Upsilon_{\bar E,\Omega^1 A}(\overline{\sigma_{\bar E}^{-1}
\nabla_{\bar E}(\bar e)})\cr
&=& (\star^{-1}\tens \id)\Upsilon_{\bar E,\Omega^1 A}(\overline{\sigma_{\bar E}^{-1}
(\star^{-1}\tens \id)\Upsilon_{E,\Omega^1 A}(\overline{\sigma_E^{-1}\nabla_E(e)})})\ ,
\end{eqnarray*}
and this should be equal to $(\id\tens \mathrm{bb}_E)\nabla_E e
=(\id\tens \mathrm{bb}_E)\sigma_E\sigma_E^{-1}\nabla_E e$. For
this to be true in general, we put $\sigma_E^{-1}\nabla_E e=f\tens\xi$, and then
we need
\begin{eqnarray*}
(\id\tens \mathrm{bb}_E)\sigma_E(f\tens\xi) &=&
(\star^{-1}\tens \id)\Upsilon_{\bar E,\Omega^1 A}(\overline{\sigma_{\bar E}^{-1}
(\star^{-1}\tens \id)\Upsilon_{E,\Omega^1 A}(\overline{f\tens\xi})})\ ,
\end{eqnarray*}
and this can be written as
\begin{eqnarray*}
\sigma_E(f\tens\xi) &=&
(\star^{-1}\tens \mathrm{bb}_E^{-1})\Upsilon_{\bar E,\Omega^1 A}(\overline{\sigma_{\bar E}^{-1}
(\xi^*\tens \bar f})\ .
\end{eqnarray*}
If we substitute into this the expression
\begin{eqnarray*}
\sigma_{\bar E}^{-1} \,=\, (\id\tens \star^{-1})\Upsilon_{\Omega^1 A,E}
\overline{\sigma_E}\Upsilon_{E,\Omega^1 A}^{-1}(*\tens \id)\ ,
\end{eqnarray*}
we the condition
\begin{eqnarray*}
\sigma_E(f\tens\xi) &=&
(\star^{-1}\tens \mathrm{bb}_E^{-1})\Upsilon_{\bar E,\Omega^1 A}
(\overline{(\id\tens \star^{-1})\Upsilon_{\Omega^1 A,E}
\overline{\sigma_E(f\tens\xi)}})\ .
\end{eqnarray*}
Substituting $\sigma_E(f\tens\xi)=\eta\tens h$, we require
\begin{eqnarray*}
\eta\tens h &=&
(\star^{-1}\tens \mathrm{bb}_E^{-1})\Upsilon_{\bar E,\Omega^1 A}
(\overline{(\id\tens \star^{-1})\Upsilon_{\Omega^1 A,E}
\overline{\eta\tens h}})  \cr
 &=&
(\star^{-1}\tens \mathrm{bb}_E^{-1})\Upsilon_{\bar E,\Omega^1 A}
\overline{(\bar h\tens\eta^*)}  \,=\,(\star^{-1}\tens \mathrm{bb}_E^{-1})
(\overline{\eta^*}\tens \overline{\bar h})\ .\quad\square
\end{eqnarray*}

\section{Riemannian geometry}

\subsection{Hermitian structures on bimodules}

There are two sensible definitions of Hermitian structures on bimodules, depending on which side
you put the bar. We have chosen this one, as it fits better with left covariant derivatives.

\begin{defin}\label{hermdeff}
A (non degenerate) Hermitian structure on an $A$-bimodule $E$ is given by an invertible morphism
$G:\overline{E}\to E^\circ$. From this we define an inner product
$\<,\>=\ev_E(\id\tens G):E\tens_A \overline{E}\to A$, and this is required to satisfy the condition
that the following composition is just $\<,\>$:
\begin{eqnarray*}
E \tens_A \overline{E} \stackrel{\mathrm{bb}\tens\id} \longrightarrow
\overline{ \overline{E}} \tens_A \overline{E} \stackrel{\Upsilon^{-1}} \longrightarrow
\overline{E \tens_A \overline{E}} \stackrel{\overline{\<,\>}} \longrightarrow \bar A
\stackrel{*^{-1}} \longrightarrow  A
\end{eqnarray*}
\end{defin}

We write $\<e,\bar f\>=\ev(e \tens_A G(\bar f))$.
We shall take a moment to see what this definition actually means.
For $e,f\in E$ the composition in \ref{hermdeff}  is
\begin{eqnarray*}
e\tens \bar f \longmapsto  \overline{\bar e} \tens \bar f
\longmapsto \overline{f\tens\bar e} \longmapsto
 \overline{\< f,\bar e\>} \longmapsto \< f,\bar e\>^*\ .
\end{eqnarray*}
Thus the condition in \ref{hermdeff} is that $\<e,\bar f\>=\<f,\bar e\>^*$. There are some other
formulae which are virtually automatic from the definition.
Since $\ev_E$ is a right $A$-module map, for all $a\in A$,
\begin{eqnarray*}
\<e,\overline{a.f}\>\,=\, \<e,\bar f.a^*\>\,=\, \<e,\bar f\>\,a^*\ .
\end{eqnarray*}
Since $\ev_E$ is a left $A$-module map, for all $a\in A$,
\begin{eqnarray*}
\<a.e,\bar f\>\,=\, a\,\<e,\bar f\>\ .
\end{eqnarray*}
Since we are using the tensor product over $A$,
\begin{eqnarray*}
\<e.a,\bar f\>\,=\, \<e,a.\bar f\>\,=\, \<e,\overline{f.a^*}\>\ .
\end{eqnarray*}

\begin{propos}
Suppose that $E$ is finitely generated projective as a left module, with dual basis
$e_i\tens e^i\in E^\circ\tens E$, and
let $G$ be a non-degenerate Hermitian structure on $E$.
Suppose that we set $g^{ij}=\<e^i,\overline{e^j}\>$, so it is automatic that
$g^{ij*}=g^{ji}$. Then we have $G(\overline{e^i})=e_j.g^{ji}$
(summation convention applies). We define
$G^{-1}(e_i)=\overline{g_{ij}.e^j}$, where without
loss of generality
we can assume that $g_{ij}.\ev(e^j\tens e_k)= g_{ik}$. Then:

a)\quad $g^{ij}\,g_{jk}\,=\, \ev(e^i\tens e_k)$\ .

b)\quad $g_{ij}\,g^{jk}\,=\, \ev(e^k \tens e_i)^*$\ .

c)\quad $g_{iq}^*\,=\, g_{qi}$\ .

\end{propos}\label{cvgads}
\proof Begin with
\begin{eqnarray*}
\overline{e^i} \,=\, G^{-1}(G(\overline{e^i}))\,=\, G^{-1}(e_j.g^{ji}) \,=\, G^{-1}(e_j) .g^{ji}\,=\,
\overline{g_{jk}.e^k} .g^{ji}  \,=\,
\overline{g^{ij}\,g_{jk}.e^k}    \  ,
\end{eqnarray*}
and apply $e_n$ to both sides to get (a). We also have
\begin{eqnarray*}
e_i \,=\, G(G^{-1}(e_i))\,=\, G(\overline{g_{ij}.e^j})\,=\,
G(\overline{e^j}).g_{ij}^*\,=\, e_k.g^{kj}\, g_{ij}^*\ ,
\end{eqnarray*}
and applying both sides to $e^p$ gives $\ev(e^k \tens e_i)=g^{kj}\, g_{ij}^*$,
 and applying $*$ to this gives
(b). Now consider
\begin{eqnarray*}
g_{ni}\,=\,
g_{nk}\,\ev(e^k \tens e_i) \,=\, g_{nk}\,g^{kj}\, g_{ij}^*\,=\,\ev(e^j \tens e_n)^*\, g_{ij}^*
\,=\,(g_{ij}\, \ev(e^j \tens e_n))^*\,=\,g_{in}^*\ .\quad\square
\end{eqnarray*}

\medskip Note that the fact that we have not defined $g_{nk}$ as the inverse to the matrix
$g^{nk}$ is nothing to do with noncommutativity. Even in ordinary differential geometry, this identification with the inverse requires choosing a chart which trivialises the bundle.

It will be convenient to define the following matrices, so that we can use matrix multiplication rather than indices:
\begin{eqnarray}\label{matrixform}
(P)_{ij}\,=\,\ev(e^i\tens e_j)\ ,\quad (g_\bullet)_{ij}=g_{ij}\ ,\quad (g^\bullet)_{ij}=g^{ij}\ .
\end{eqnarray}
Then, just from the definition of finitely generated projective, $P^2=P$. Then the results of
\ref{cvgads} can be summarised as
\begin{eqnarray}\label{matrixform77}
g^{\bullet*}=g^\bullet \ ,\quad g_\bullet^*=g_\bullet\ ,\quad
g^\bullet g_\bullet=P\ ,\quad g_\bullet P=g_\bullet\ ,\quad Pg^\bullet= g^\bullet\ .
\end{eqnarray}

\begin{propos} For all $a\in A$, if we set the matrix $(P(a))_{ij}=\ev(e^i.a\tens e_j)$, then
$g^\bullet\,P(a)^*=P(a^*)\ g^\bullet$.
\end{propos}
\proof
\begin{eqnarray*}
g^{ij}\, \ev(e^k.a\tens e_j)^* &=& \<e^i,\overline{e^j}\>\, \ev(e^k.a\tens e_j)^* \cr
&=& \<e^i,\overline{\ev(e^k.a\tens e_j)\,e^j}\> \cr
&=& \<e^i,\overline{e^k.a}\>       \cr
&=&  \<e^i,a^*.\overline{e^k}\> \cr
&=&  \<e^i.a^*,\overline{e^k}\>  \cr
&=&   \<\ev(e^i.a^*\tens e_j).e^j,\overline{e^k}\>     \cr
&=& \ev(e^i.a^*\tens e_j)\ \<e^j,\overline{e^k}\> \ .\quad\square
\end{eqnarray*}

\subsection{Christoffel symbols}
Begin with a left covariant derivative $\nabla$ on a  right $A$-module $E$. We suppose that $E$ is
finitely generated projective as a left $A$-module, with dual basis $e^i\in E$
and $e_i\in E^\circ$. Then we define the Christoffel symbols
\begin{eqnarray*}
\Gamma^j_i \,=\,-\, (\id\tens\ev)( \nabla e^j\tens e_i)\in\Omega^1 A\ .
\end{eqnarray*}
(We choose
the minus sign to fit with the standard convention for the covariant derivative
of 1-forms, and the reader should remember that the basis
of the 1-forms should be written with upper indices if the coefficients
of a 1-form have lower indices, as is standard.)
We make the Christoffel symbols into a matrix by defining
\begin{eqnarray*}
(\Gamma)_{ji}\,=\, \Gamma^j_i \ .
\end{eqnarray*}

\begin{propos}\label{vchuwb}
The matrix Christoffel symbols obey the following equations:
\begin{eqnarray*}
\Gamma\,P\,=\,\Gamma\ ,\quad \Gamma\,=\,P\,\Gamma-\extd P.P\ .
\end{eqnarray*}
\end{propos}
\proof From the definition of the Christoffel symbols we have
$\Gamma_k^i.\ev( e^k\tens e_j)=\Gamma_j^i$.
From the definition  of finitely generated projective,
\begin{eqnarray*}
\nabla(e^i) &=& \nabla(\ev(e^i\tens e_j).e^j) \cr
&=& \ev(e^i\tens e_j).\nabla_E(e^j)+\extd( \ev(e^i\tens e_j)) \tens e^j\ .\quad\square
\end{eqnarray*}

\medskip In the same manner as the Christoffel symbols, we can describe $\sigma$
(in the cases where it exists) by the matrix
\begin{eqnarray*}
\Sigma_{ij}(a) &=& (\id\tens\ev)(\sigma(e^i\tens\extd a)\tens e_j)\in\Omega^1 A\ .
\end{eqnarray*}

\begin{propos}
If $(\nabla,E,\sigma)$ gives a bimodule connection, then
\begin{eqnarray*}
\Sigma(a) &=& \extd P(a).P-P(a)\,\Gamma+\Gamma\,P(a)\ .
\end{eqnarray*}
\end{propos}
\proof
\begin{eqnarray*}
\nabla(e^i.a) &=& \nabla(\ev(e^i.a\tens e_j).e^j) \cr
&=&  \extd\ \ev(e^i.a\tens e_j)\tens e^j  + \ev(e^i.a\tens e_j).\nabla(e^j) \ ,
\end{eqnarray*}
so we have
\begin{eqnarray*}
(\id\tens\ev)(\nabla(e^i.a)\tens e_k-\nabla(e^i).a\tens e_k) &=&
\extd\ \ev(e^i.a\tens e_j) . \ev(e^j\tens e_k)-\ev(e^i.a\tens e_j).\Gamma^j_k \cr
&& +\, \Gamma^i_j.\ev(e^j.a\tens e_k)\ .\quad\square
\end{eqnarray*}

\subsection{Connections preserving Hermitian structures}
Using the map $G^{-1}:E^\circ\to \bar E$
we have an element $(G^{-1}\tens\id)\coev^R_E(1)\in \overline{E}\tens E$
which commutes with elements of $A$ (as $G^{-1}$ is a bimodule map),
and we wish to find the conditions for it to have zero covariant derivative.

\begin{propos} \label{bcadhsucv}
\begin{eqnarray*}
&& ( \sigma_{\overline{E}}^{-1}  \tens\id)\nabla_{\bar E\tens E}\big((G^{-1}\tens\id)\coev^R_E(1)\big)\cr
&=&
-\,\overline{e^k}\tens (\Gamma^j_k)^* .g_{ji}\tens e^i+
 \overline{e^j} \tens\extd g_{ji}\tens e^i- \overline{e^j} .g_{ji}\tens\Gamma_k^i\tens e^k
 \in \overline{E}\tens_A \Omega^1 A\tens_A E\ .
\end{eqnarray*}
\end{propos}
\proof First
\begin{eqnarray*}
(G^{-1}\tens\id)\coev^R_E(1) \,=\, (G^{-1}\tens \id)(e_i\tens e^i)\,=\,
\overline{g_{ij}.e^j} \tens e^i\,=\, \overline{e^j} .g_{ji}\tens e^i\ ,
\end{eqnarray*}
and on applying $\nabla$ to this we get
\begin{eqnarray*}
\nabla_{\overline{E}} (\overline{e^j}) .g_{ji}\tens e^i+(\sigma_{\overline{E}} \tens\id)
( \overline{e^j} \tens\extd g_{ji}\tens e^i+ \overline{e^j} .g_{ji}\tens \nabla_E(e^i))\ ,
\end{eqnarray*}
and applying $ \sigma_{\overline{E}}^{-1}  \tens\id $ to this gives
\begin{eqnarray} \label{fgusat}
 \sigma_{\overline{E}}^{-1} \nabla_{\overline{E}} (\overline{e^j}) .g_{ji}\tens e^i+
 \overline{e^j} \tens\extd g_{ji}\tens e^i+ \overline{e^j} .g_{ji}\tens \nabla_E(e^i)\ .
\end{eqnarray}
Now use the description of $\nabla_{\overline{E}}$ and $ \sigma_{\overline{E}}$ in  \ref{sliogcb}
to show that (\ref{fgusat}) can be rewritten as
\begin{eqnarray} \label{fgusarr}
(\id\tens\star^{-1})\,\Upsilon\, \overline{\nabla_{E} (e^j)} .g_{ji}\tens e^i+
 \overline{e^j} \tens\extd g_{ji}\tens e^i+ \overline{e^j} .g_{ji}\tens \nabla_E(e^i)\ .
\end{eqnarray}
Now use
\begin{eqnarray*}
\nabla_E e^i \,=\, -\,\Gamma_k^i\tens e^k
\end{eqnarray*}
to rewrite (\ref{fgusarr}) as
\begin{eqnarray*}
-\,\overline{e^k}\tens (\Gamma^j_k)^* .g_{ji}\tens e^i+
 \overline{e^j} \tens\extd g_{ji}\tens e^i- \overline{e^j} .g_{ji}\tens\Gamma_k^i\tens e^k\ .\quad\square
\end{eqnarray*}

\begin{propos}\label{ksdzjcv}
The condition for a connection to preserve the Hermitian metric is
\begin{eqnarray*}
g_\bullet.\Gamma&=&
\frac12 P^*.\extd g_\bullet.P+\phi\ ,
\end{eqnarray*}
where $\phi\in M_n(\Omega^1 A)$ with $\phi^*=-\phi$,
$P^*\phi=\phi$ and  $\phi\,P=\phi$. From this we can deduce that
\begin{eqnarray*}
\Gamma&=&
\frac12 \,g^\bullet.\extd g_\bullet.P+g^\bullet.\phi-\extd P.P\ .
\end{eqnarray*}
\end{propos}
\proof Supposing that $\nabla_{\bar E\tens E}\big((G^{-1}\tens\id)\coev^R_E(1)\big)
=0$, taking the inner product of both sides of the result of \ref{bcadhsucv}
\begin{eqnarray*}
g_\bullet\, g^\bullet.\Gamma^*.P^*\, g_\bullet+g_\bullet \,P.\Gamma.g^\bullet \, g_\bullet&=&
g_\bullet\, g^\bullet.\extd g_\bullet.g^\bullet\, g_\bullet\ ,
\end{eqnarray*}
and using (\ref{matrixform77}) gives
\begin{eqnarray*}
P^*.\Gamma^*. g_\bullet+g_\bullet.\Gamma.P&=&
P^*.\extd g_\bullet.P\ ,
\end{eqnarray*}
and then using \ref{vchuwb} gives
\begin{eqnarray*}
\Gamma^*. g_\bullet+g_\bullet.\Gamma&=&
P^*.\extd g_\bullet.P\ .
\end{eqnarray*}
From this
\begin{eqnarray*}
g_\bullet.\Gamma&=&
\frac12 P^*.\extd g_\bullet.P+\phi\ ,
\end{eqnarray*}
where $\phi\in M_n(\Omega^1 A)$ with $\phi^*=-\phi$. From (\ref{matrixform77}) again we see that
$P^*\,\phi=\phi$, and thus that $\phi\,P=\phi$. Multiplying this on the left by $g^\bullet$ gives
\begin{eqnarray*}
P.\Gamma&=&
\frac12 g^\bullet.\extd g_\bullet.P+g^\bullet.\phi\ ,
\end{eqnarray*}
and \ref{vchuwb} gives the answer.\quad$\square$

\medskip Now we look at a weaker condition than full metric preservation introduced in
\cite{Ma:rief} -- that of vanishing cotorsion. As in the case of torsion, this only applies in the case
$E=\Omega^1 A$. We modify it slightly here to fit in with the bar notation for left connections
(originally it took values in $\Omega^2 A \tens_A \Omega^1 A$).

\begin{defin}\label{kxzcvzvck}
If $\nabla$ is a left covariant derivative on $\Omega^1 A$
 and $G$ is a Riemannian metric on $\Omega^1 A$, define the cotorsion as
\begin{eqnarray*}
(\id\tens\wedge)( \sigma_{\overline{\Omega^1}}^{-1}  \tens\id)\nabla_{\overline{ \Omega^1}\tens \Omega^1}\big((G^{-1}\tens\id)\coev^R_{\Omega^1}(1)\big)\in\overline{\Omega^1 A}\tens_A\Omega^2 A\ .
\end{eqnarray*}
\end{defin}

Just to check that this is a weaker condition than full metric preservation,
the reader should note that if the metric is preserved by $\nabla$, then the formula in
\ref{kxzcvzvck} becomes
\begin{eqnarray*}
(\id\tens\wedge)( \sigma_{\overline{\Omega^1}}^{-1}  \tens\id)
\,(\id\tens G^{-1}\tens\id)
\nabla_{( \Omega^1)^\circ\tens \Omega^1}\big(\coev^R_{\Omega^1}(1)\big)\in\overline{\Omega^1 A}\tens_A\Omega^2 A\ ,
\end{eqnarray*}
and this vanishes by definition of $\nabla_{( \Omega^1)^\circ}$.

\begin{propos}\label{cvcvictcttcttc}
In terms of a dual basis and using matrix notation, the condition for vanishing cotorsion is
that the following expression vanishes:
\begin{eqnarray*}
P^*\,\extd g_\bullet\wedge e^\bullet-g_\bullet (\extd P.P+\Gamma)\wedge e^\bullet-\Gamma^* g_\bullet
\wedge e^\bullet\ .
\end{eqnarray*}
\end{propos}
\proof By \ref{bcadhsucv} the cotorsion is
\begin{eqnarray*}
-\,\overline{e^j}\tens (\Gamma^p_j)^* .g_{pi}\wedge e^i+
 \overline{e^j} \tens\extd g_{ji}\wedge e^i- \overline{e^j} .g_{ji}\tens\Gamma_k^i\wedge e^k\ .
\end{eqnarray*}
Applying $\<e^s,-\>$ to this gives
\begin{eqnarray*}
-\,g^{sj}\, (\Gamma^p_j)^* .g_{pi}\wedge e^i+
 g^{sj}\, \extd g_{ji}\wedge e^i- g^{sj}\, g_{ji}\,\Gamma_k^i\wedge e^k\ .
\end{eqnarray*}
Multiplying on the left by $g_{rs}$ and using (\ref{matrixform77})  gives
\begin{eqnarray*}
P^*\,\extd g_\bullet\wedge e^\bullet-g_\bullet P\,\Gamma\wedge e^\bullet-P^*\Gamma^* g_\bullet
\wedge e^\bullet\ .
\end{eqnarray*}
Now from \ref{vchuwb}, $P\,\Gamma=\extd P.P+\Gamma$ and $\Gamma\,P=\Gamma$, so
we get the answer.\quad$\square$

\medskip The reader may check that the matrix expression in \ref{cvcvictcttcttc} vanishes,
given the matrix condition for metric compatibility in \ref{ksdzjcv}.

\subsection{Star compatible connections}
Now suppose that $E$ is a star object in our bar category ${}_A\CM_A$ in the sense of
 \cite{barcats}, i.e.\ we require a morphism $\star:E\to \bar E$ with $\bar\star\star(e)=\bar{\bar e}$ for all $e\in E$. If $\nabla$ is a bimodule covariant derivative we can require that it preserves the star map in the sense
\begin{eqnarray}\label{bcihdbbdiuv}
(\id\tens\star)\nabla_{E}\,=\,\nabla_{\overline{E}}\,\star:E\to
\Omega^1A \tens_{A} \overline{E}\ .
\end{eqnarray}
that $\star$ is a morphism in ${}_A\mathcal{E}_A$. This is equivalent to saying that $\star$ extends to make $(E,\nabla,\sigma)$ a star-object.  In fact it will be convenient to consider the potentially weaker morphism condition in  \ref{bchdjsabv}, which we call `star-compatible':

\begin{defin}\label{badsuvcvbi} If $E$ is a star-object in ${}_A\CM_A$, we say that $\nabla$ is {\em star compatible} if
\begin{eqnarray*}
(\id\tens\star)\sigma_{E}=\sigma_{\overline{E}}(\star\tens\id):
E \tens_A\Omega^1 A \to \Omega^1 A \tens_A \overline{E}\ .
\end{eqnarray*}
\end{defin}

\medskip We will use this idea later in the paper, rather that that of $\star$ being a morphism in the category ${}_A\mathcal{E}_A$.

\begin{propos}
Using theorem \ref{sliogcb}, we can rephrase the condition for star compatibility as
\begin{eqnarray*}
\sigma_{E}&=&  (\star^{-1}\tens\star^{-1}) \Upsilon\, \overline{\sigma_E^{-1}}\,
\Upsilon^{-1} (\star\tens\star)\ .
\end{eqnarray*}
The condition for $\star$ to be a morphism in the category ${}_A\mathcal{E}_A$ is that
\begin{eqnarray*}
\overline{\sigma_E}\,\Upsilon^{-1}\,(\star\tens\star)\,\nabla_E(e) &=& \overline{\nabla_E(e^*)}\ .
\end{eqnarray*}
\end{propos}
\proof From theorem \ref{sliogcb} we have
\begin{eqnarray*}
\sigma_{\overline{E}}(\star\tens\id) &=&  (\star^{-1}\tens\id) \Upsilon\, \overline{\sigma_E^{-1}}\,
\Upsilon^{-1} (\star\tens\star)\ ,
\end{eqnarray*}
so the condition that the connection is star compatible is that
\begin{eqnarray*}
\sigma_{E}&=&  (\star^{-1}\tens\star^{-1}) \Upsilon\, \overline{\sigma_E^{-1}}\,
\Upsilon^{-1} (\star\tens\star)\ .
\end{eqnarray*}
Also from theorem \ref{sliogcb} we have
\begin{eqnarray*}
(\id\tens\star^{-1})\nabla_{\overline{E}}\,\star(e) &=& (\id\tens\star^{-1})
\nabla_{\overline{E}}(\overline{e^*}) \cr
&=& (\star^{-1}\tens\star^{-1}) \Upsilon\,\overline{\sigma_E^{-1}
\nabla_E(e^*)}\ .\quad\square
\end{eqnarray*}

\medskip We shall later be considering the case $E=\Omega^1 A$, in which case we shall use the
star operation discussed at the beginning of \ref{ikcvsaucfv}.

\section{Bimodule covariant derivatives on Hopf *-algebras}

In this section we specialise to the case where our `coordinate algebra' is a `quantum group' or Hopf algebra $H$. This has a coproduct $\Delta:H\to H\tens H$, a counit $\eps:H\to \C$ and an antipode $S:H\to H$ with the usual axioms expressing a `group' structure on the coordinate algebra. We refer to
\cite{Ma:prim} for details. We use the shorthand notation $\Delta a=a\o\tens a\t$ for output of $\Delta$. A Hopf algebra coacts on itself from both the left and the right via the coproduct, corresponding to left and right translation in the `group'.

We start with a preliminary subsection recalling basic facts about the differential calculi on the quantum group $\C_q[SU_2]$ and functions on a finite group, that will be needed in our classification results. We then study bimodule covariant derivatives in terms of the left invariant part $\nabla^L$ of the connection, which is the new part of the present section. Section~4 then provides detailed results for our chosen examples.

\subsection{Preliminary background}
We suppose that $\Omega^1H$ is left-covariant in the sense that  there is a well defined map
\[ \Delta_L:\Omega^1 H \to H\tens \Omega^1 H,\quad a.\extd b \mapsto a_{(1)}\, b_{(1)}\tens a_{(2)}.\extd b_{(2)},\]
 which we will write as $\Delta_L\xi=\xi\mapsto \xi_{[-1]}\tens \xi_{[0]}$. In this case $\Delta_L$
is a left $H$-coaction and a bimodule map. As with the theory of Lie groups, we can then
  trivialise the module of 1-forms by the linear isomorphism
$Y:\Omega^1 H \to H\tens \Lambda^1 H$ defined by
$Y(\xi)=\xi_{[-2]}\tens S(\xi_{[-1]}) \xi_{[0]}$, where
$\Lambda^1 H$ is the space of left invariant 1-forms. The bimodule structure of $\Omega^1H$ appears now as the right and left module structures on $H\tens \Lambda^1 H$ given by
\begin{eqnarray*}
h.(g\tens \xi) \,=\, h\,g\tens \xi\ ,\quad (g\tens \xi).h \,=\, g\, h_{(1)}\tens \xi\ra h_{(2)}\ ,
\end{eqnarray*}
where the right action on $\Lambda^1 H$ is given by $\xi\ra h =
S(h_{(1)}).\xi.h_{(1)}$. We can also specify the left invariant derivative
$\varpi:H\to \Lambda^1 H$ by $\varpi(h)=S(h_{(1)}).\extd h_{(2)}$. It is conventional to restrict this map to $H^+=\ker\eps$ (the augmentation ideal) since $H=\C1\oplus H^+$ and $\varpi$ is trivial on $1$. In this case
\[ \Lambda^1 H\ \isom\ H^+/I_H,\quad I_H=\ker\varpi\]
where $I_H$ is a right ideal. Left covariant $\Omega^1H$ are in 1-1-correspondence with right ideals of $H^+$ in this way. The  right action $\ra$ on $\Lambda^1H$ in this form is inherited from the product of $H$. In the finite codimension case the vector space
 dual  $(\Lambda^1H)'=\mathfrak{h}$ is called the `quantum Lie algebra' associated to the differential calculus on $H$. It can be viewed in this way as the space of left-invariant vector fields as a subset of ${\rm Vec}(H)$.

\medskip If there is also a well defined map
$\Delta_R:\Omega^1 H \to  \Omega^1 H\tens H$ given by
$a.\extd b \mapsto a_{(1)}.\extd b_{(1)}\tens a_{(2)}\, b_{(2)}$, we call the
differential calculus bicovariant. (This just means that the coproduct is
differentiable.) Where it exists, it is a right coaction, and we denote it by
$\Delta_R(\xi)=\xi_{[0]}\tens \xi_{[1]}$.
  Of particular importance in this case is that the right action and right adjoint coaction make $\Lambda^1 H$ a `crossed' or `Drinfeld-Radford-Yetter' module (i.e. a module under the quantum double of $H$) and as such there is an induced braiding on it,
\begin{eqnarray*}
\Psi(\xi\tens\eta) &=& \eta_{[0]}\tens \xi\ra\eta_{[1]}
\end{eqnarray*}
which obeys the braid relations and is invertible if the antipode $S$ is. We refer to \cite{Ma:prim} for an introduction and to \cite{Wor:dif} for the original work.  In terms of the ideal $I_H$, the calculus is bicovariant if and only if $I_H$ is stable under the right adjoint coaction $\Ad(a)=a\t\tens (Sa\o)a\thr$. The coaction $\Delta_R$ on $\Lambda^1H$ in the form  $H^+/I_H$ is inherited from this. Hence one can also write
\[ \Psi(\xi\tens\varpi(a))=\varpi(a\o)\tens\xi\ra (Sa\o)a\t,\quad \forall a\in H.\]
Note that $\Psi$ might be well defined even if the calculus is not bicovariant. Moreover if $\Psi$ is well defined, it implies a well defined map $\Omega^1 H\tens_H\Omega^1 H \to \Omega^1 H\tens_H\Omega^1 H $ obeying the braid relations. We define $\Omega^2 H$ as the skew symmetrisation of
$\Omega^1 H\tens_H\Omega^1 H$ with respect to any such $\Psi$, more precisely
$\Omega^2$ is $\Omega^1 H\tens_H\Omega^1 H$ divided by the kernel of $\id-\Psi$.
Although we will usually talk about 1-forms, we will assume that 2-forms are defined in this way.

Finally, all our constructions are consistent with $H$ being a Hopf $*$-algebra. Here $H$ is a $*$-algebra, $\Delta$ is a $*$-algebra map and $(S\circ *)^2=\id$. In this case $\Omega^1H$ is a $*$-object as in the previous section if and only if $I_H$ is stable under $*$. In this the coactions are
automatically compatible under the star structure.

\subsection{Left invariant covariant derivatives}
Using the left Liebniz rule
 we can
reduce specifying a covariant derivative
 $\nabla:\Omega^1 H \to \Omega^1 H\tens_H \Omega^1 H$
 to specifying it
on the left invariant forms:
\begin{eqnarray} \label{reconner}
\nabla(\eta) &=& \nabla(\eta_{[-2]}S(\eta_{[-1]}).\extd\eta_{[0]}) \cr
&=& \extd \eta_{[-2]} \tens_H S(\eta_{[-1]}).\extd\eta_{[0]} +
\eta_{[-2]}.\nabla^L(S(\eta_{[-1]}).\extd\eta_{[0]})\ .
\end{eqnarray}
Here we write the restriction of the covariant derivative
as $\nabla^L:\Lambda^1 H \to \Omega^1 H\tens_H \Omega^1 H$.
The map $\sigma$ associated to the covariant derivative, if it exists,
is necessarily given by the following formula in terms of
 $\nabla^L$,
where $\xi\in \Lambda^1 H$ and $h\in H$:
\begin{eqnarray} \label{jzxdhdhvbc}
\sigma(\xi\tens\extd h_{(2)}.S^{-1}(h_{(1)})) &=&
\sigma(\xi\tens\extd h_{(2)}).S^{-1}(h_{(1)}) \cr
&=& (\nabla(\xi. h_{(2)})-\nabla(\xi).h_{(2)}).S^{-1}(h_{(1)}) \cr
&=& \nabla(\xi. h_{(2)}).S^{-1}(h_{(1)})-\nabla(\xi).\epsilon(h) \cr
&=& \nabla(h_{(2)}S(h_{(3)}).\xi. h_{(4)}).S^{-1}(h_{(1)})-\nabla(\xi).\epsilon(h) \cr
&=& \nabla(h_{(2)}(\xi\ra h_{(3)})).S^{-1}(h_{(1)})-\nabla(\xi).\epsilon(h) \cr
&=& \extd h_{(2)} \tens (\xi\ra h_{(3)}).S^{-1}(h_{(1)})
-\nabla(\xi).\epsilon(h) \cr && +\,
h_{(2)} .\nabla(\xi\ra h_{(3)}).S^{-1}(h_{(1)}) \cr
&=& \extd h_{(3)}.S^{-1}(h_{(2)}) \tens \xi\ra h_{(4)}S^{-1}(h_{(1)})
-\nabla^L(\xi).\epsilon(h) \cr && +\,
\nabla^L(\xi\ra h_{(2)})\ra S^{-1}(h_{(1)}) \ .
\end{eqnarray}
It is rather neater to restate this in terms of $\varpi$:
\begin{eqnarray} \label{varpiytyty}
\sigma(\xi\tens\varpi(h)) &=&
\varpi(h_{(2)}) \tens \xi\ra S(h_{(1)})\, h_{(3)}
+\nabla^L(\xi).\epsilon(h)-\nabla^L(\xi\ra S(h_{(1)}) )\ra h_{(2)}\ ,
\end{eqnarray}
The problem here is whether $\sigma$ is a well defined
function of its second variable. We recognise the first term $\Psi$ here as something that is
well-defined at least in the bicovariant case and some other cases. Thus we shall write
\begin{eqnarray} \label{varpiytyty77}
\sigma&=& \Psi+\hat\sigma\ ,
\end{eqnarray}
and we note that $\sigma$ is well defined if and only $\hat\sigma$ is well defined.

Now we can ask for $\nabla$ to be left covariant. This is the same as asking if
a covariant derivative on a Lie group is invariant to left translation by group
elements. This is defined by the commutativity of the following diagram,
where the dots indicate that we use the tensor product left coaction:

\begin{picture}(100,80)(-100,6)

\put(60,63){\vector(1,0){45}}
\put(60,23){\vector(1,0){45}}
\put(28,53){\vector(0,-1){20}}
\put(137,53){\vector(0,-1){20}}
\put(17,60){$\Omega^1 H$}
\put(-27,20){${}^\bullet\Omega^1 H \tens_H {}^\bullet\Omega^1 H$}
\put(111,60){$H\tens \Omega^1 H $}
\put(111,20){$ H \tens \Omega^1 H \tens_H \Omega^1 H$}

\put(75,70){$\nabla^L$}
\put(75,30){$\nabla^L$}
\put(12,40){$\nabla$}
\put(143,40){$\id\tens\nabla$}
\end{picture}

\begin{propos}
The covariant derivative $\nabla$ is left covariant if and only if it restricts to a map
$\nabla^L:\Lambda^1 H \to \Lambda^1 H\tens \Lambda^1 H$.
Furthermore we see from (13) that if if we have a bomodule covarariant derivative $\nabla$ which is left covariant, then $\sigma$ also restricts to a map from $\Lambda^1 H\tens \Lambda^1 H$ to itself.
\end{propos}
\proof By writing $\xi\tens_H\eta=\xi.\eta_{[-2]}\tens_H S(\eta_{[-1]}).\eta_{[0]}
\in\Omega^1 H\tens_H \Omega^1 H$.\quad$\square$

\subsection{Bicovariant differential calculi}
The most obvious choice of left invariant connection is given by
$\nabla^L=0$, and this gives the left Maurer-Cartan connection.
Clearly this extends to a bimodule covariant derivative on $H$ if and only if $\sigma=\Psi$ is
well-defined.  A sufficient condition is that the calculus is bicovariant. More generally in the bicovariant case,  considering $\nabla^L:\Lambda^1 H \to \Lambda^1 H\tens \Lambda^1 H$,
we can write some immediate left covariant bimodule connections as follows:

\medskip\noindent
1)\quad $\nabla^L(\xi)=0$, the left Maurer Cartan connection,
$\sigma=\Psi$.

\medskip\noindent
2)\quad $\nabla^L(\xi)=-\varpi(S^{-1}(\xi_{[1]}))\tens \xi_{[0]}$,
 the right Maurer Cartan connection,
$\sigma=\Psi^{-1}$.

\medskip\noindent
3)\quad $\nabla^L(\xi) =\xi_{[0]}\tens \varpi(\xi_{[1]})$,
$\sigma= \Psi+\id-\Psi^2$.

\begin{remark} Consider a Hopf algebra $H$ with a bicovariant differential calculus.
If $(E,\nabla,\sigma)$ gives a bimodule connection and
if $\nabla^L$ is right covariant, then $\sigma$ is a map in the crossed module or Drinfeld-Radford-Yetter  category. This is because all
we have to check is that it preserves the right action, and it does this as it is a bimodule map.
This is the case in the above three examples.
\end{remark}

\section{Example: The function algebra of a finite group}

For $G$ a finite group, we use the usual Hopf algebra $C(G)$, with basis
$\delta_x$ for $x\in G$
(the function taking value 1 at $x$
and zero elsewhere). This has operations
\begin{eqnarray*}
\delta_x.\delta_y\,=\, \delta_{x,y}\, \delta_x\ ,\quad
\Delta\delta_x\,=\,\!\!\!\sum_{y,z\in G:\,yz=x}\!\!\delta_y\tens\delta_z\ ,\quad
1\,=\, \sum_{x\in G}\delta_x\ ,\quad \epsilon(\delta_x)\,=\,\delta_{x,e}\ ,\quad
S(\delta_x)\,=\, \delta_{x^{-1}}.
\end{eqnarray*}
Here $e\in G$ is the identity element, and $\delta_{x,y}$ is the usual
Kroneker delta.
For $f\in C(G)$ and $g\in G$, it will be convenient to define
the right translation $R_g(f)\in C(G)$
by $R_g(f)(x)=f(xg)$ so that $R_g(\delta_x)=\delta_{xg^{-1}}$. The star operation on $C(G)$ is given by
\begin{eqnarray*}
\star\delta_x \,=\, \overline{\delta_x}\ .
\end{eqnarray*}

We give $C(G)$ a differential calculus as follows \cite{Ma:rief}:
Take $C$ to be a subset of $G$
which does not include the identity. Then take the left invariant
1-forms to have basis $\xi^c$ for $c\in C$. The bimodule commutation relations and the
exterior derivative  are
\begin{eqnarray} \label{kjhcvzvahjks}
\xi^c.f\,=\, (R_c f).\xi^c\ ,\quad
\extd f \,=\,\sum_{c\in C}(R_c f-f).\xi^c\ .
\end{eqnarray}
We can invert this to give
\begin{eqnarray} \label{bhvdsbv}
\xi^c &=& \sum_{u\in G}\delta_{uc^{-1}}.\extd\delta_u\ .
\end{eqnarray}
The calculus is bicovariant if and only if $C$ is ad-stable. The right action and (in the bicovariant case) the right coaction and induced braiding are given by
\begin{equation}\label{YDhfhdjn}
 \xi^a \ra \delta_g = \delta_{a,g}\,\xi^a\, \quad \Delta_R\xi^c
=  \sum_{y\in G} \xi^{ycy^{-1}}\tens
\delta_y \ ,\quad
\Psi(\xi^a\tens\xi^b) =
 \xi^{aba^{-1}}\tens \xi^a\ .
\end{equation}
We also have
\begin{eqnarray*}
\varpi(\delta_g) &=& \sum_{xy=g}S(\delta_x)\, \extd(\delta_y)\,=\, \sum_{c\in C}  (\delta_{g,c}-\delta_{g,e}).\xi^c
\end{eqnarray*}
Thus $\varpi$ has kernel with basis $\delta_g$ for $g\in G\setminus(C\cup\{e\})$
and $\sum_{c\in C}\delta_c +\delta_e$.

If $C$ is closed under inverse, we define $\xi^{a*}=-\xi^{a^{-1}}$. Then we have
$(\extd f)^*=\extd f^*$, as:
\[
(\extd \delta_x)^* = \sum_{c\in C}((\delta_{xc^{-1}}-\delta_x).\xi^c)^* = -\sum_{c\in C} \xi^{c^{-1}}.(\delta_{xc^{-1}}-\delta_x) = \sum_{c\in C} (\delta_{xc}-\delta_{x}).\xi^{c^{-1}} = \extd\delta_x \,=\, \extd\delta_x^*\ .
\]

A basis of $(\Lambda^1 C(G))^\circ$ is given by
$\xi_c$ for $c\in C$, where we define
$\ev(\xi^a \tens \xi_c)=\delta_{c,a}$.
 The action and coaction are given by standard results on the dual
of crossed or Drinfeld-Radford-Yetter modules as
\begin{eqnarray*}
\xi_a\ra \delta_b\,=\,\delta_{a,b^{-1}}\,\xi_a\ ,\quad \Delta_R(\xi_a) \,=\,
\sum_g \xi_{gag^{-1}}\tens\delta_g\ .
\end{eqnarray*}
Moreover the calculus is inner with respect to $\theta=\sum_{a\in C} \xi^a$ in the sense that
$\extd$ is given by a graded commutator
$d=[\theta,-\}$.
Then the exterior derivative on 1-forms is given by
\begin{eqnarray*}
\extd \xi^c &=& \sum_{a\in C} \Big(\xi^a\wedge\xi^c+\xi^c\wedge\xi^a\Big)-\sum_{b,a\in C}
\delta_{c,ab}\,\xi^a\wedge\xi^b\ .
\end{eqnarray*}

\begin{propos}
A left invariant Hermitian structure can be written as $G:\overline{\Lambda^1 C(G)}\to
(\Lambda^1 C(G))^\circ$ given by $G(\overline{\xi^a})=\xi_b.g^{b,a}$, where $g^{b,a}\in\mathbb{C}$.
Then:

\smallskip\noindent
1)\quad If $G$ is a right module map, then $g^{a,b}\neq 0$ only
 if $a=b$, i.e.\ the metric is diagonal in our
basis.

\smallskip\noindent
2)\quad If $G$ is a right comodule map, then for every $a\in C$ and $x\in G$,
$g^{xax^{-1},xax^{-1}}=g^{a,a}$.
\end{propos}
\proof For (1) we have the equality of
\begin{eqnarray*}
G(\overline{\xi^a}\ra \delta_g) &=& G(\overline{\xi^a \ra \delta_{g^{-1}}})
\,=\, \delta_{a,g^{-1}}\,G(\overline{\xi^a })  \,=\,
 g^{ab}\,\delta_{a,g^{-1}}\, \xi_b \ , \cr
G(\overline{\xi^a})\ra \delta_g &=& g^{ab}\, \xi_b \ra \delta_g  \,=\,
 g^{ab}\,\delta_{b,g^{-1}}\, \xi_b\ ,
\end{eqnarray*}
so we deduce that $g^{ab}\neq 0$ only if $a=b$.

For (2), consider the right invariance property for $G$, which gives the equality
of the following:
\begin{eqnarray*}
(G\tens\id)\rho(\overline{\xi^a}) &=& \sum_{x\in G}
G(\overline{\xi^{xax^{-1}}}) \tens\delta_{x} \cr
&=&  \sum_{x\in G,\, b\in C} g^{xax^{-1},b}\, \xi_b \tens\delta_{x} \cr
&=&  \sum_{x\in G,} g^{xax^{-1},xax^{-1}}\, \xi_{xax^{-1}} \tens\delta_{x} \ , \cr
\rho\, G(\overline{\xi^a}) &=& g^{ab}\, \rho(\xi_b) \cr
&=& \sum_{x\in G,\, b\in C} g^{ab}\, \xi_{xbx^{-1}}\, \tens\delta_x\,=\,
\sum_{x\in G,} g^{aa}\, \xi_{xax^{-1}}\, \tens\delta_x\ . \quad \square
\end{eqnarray*}

\subsection{Covariant derivatives on finite groups}
For the given differential calculus on $C(G)$,
we set $\Gamma^a_{c}=\hat\Gamma^a_{bc}\, \xi^b$. If the covariant derivative
$\nabla$ is is a left comodule map, then we see that we can take $\hat\Gamma^a_{bc}\in\mathbb{C}$, as the $\xi^b$ form a basis of the left invariant 1-forms.

\begin{propos} \label{iusvggiuav}
The left invariant covariant derivative on $C(G)$ given by
\begin{eqnarray*}
\nabla^L(\xi^a)\,=\,-\,\hat\Gamma^a_{bc}\, \xi^b\tens \xi^c\ ,
\end{eqnarray*}
is a bomodule covariant derivative
if and only if
\begin{eqnarray*}
a^{-1}bc\notin C\cup\{e\}&\Rightarrow& \hat\Gamma^a_{bc}=0\ .
\end{eqnarray*}
In this case $\sigma$ is given by (summing over $b,c\in C$)
\begin{eqnarray*}
\sigma(\xi^d\tens\xi^k) &=&  \delta_{bc,dk}\,
(\hat\Gamma^d_{bc}+ \delta_{d,c}
)\, \xi^b\tens \xi^c\ ,
\end{eqnarray*}
\end{propos}
\proof
If we set $\sigma=\Psi
+\hat\sigma$, then
\begin{eqnarray*}
\hat\sigma(\xi^a\tens\varpi(\delta_g)) &=&
\nabla^L(\xi^a).\epsilon(\delta_g)-
\nabla^L(\xi^a\ra S((\delta_g)_{(1)}) )\ra (\delta_g)_{(2)} \cr
&=& \nabla^L(\xi^a)\,\delta_{g,e}-
\sum_{xy=g}\nabla^L(\xi^a\ra \delta_{x^{-1}} )\ra \delta_y \cr
&=& \nabla^L(\xi^a)\,\delta_{g,e}-
\sum_{xy=g}\delta_{a,x^{-1}}\,\nabla^L(\xi^a )\ra \delta_y \cr
&=& \nabla^L(\xi^a)\,\delta_{g,e}-\nabla^L(\xi^a )\ra \delta_{ag} \cr
&=& -\,\hat\Gamma^a_{bc}\, \xi^b\tens \xi^c\,\delta_{g,e}
+\hat\Gamma^a_{bc}\, (\xi^b\tens \xi^c)\ra \delta_{ag} \cr
&=& -\,\hat\Gamma^a_{bc}\, \xi^b\tens \xi^c\,\delta_{g,e}
+\hat\Gamma^a_{bc}\, \sum_{xy=ag}\xi^b\ra \delta_x\tens \xi^c\ra \delta_y \cr
&=& -\,\hat\Gamma^a_{bc}\, \xi^b\tens \xi^c\,\delta_{g,e}
+\hat\Gamma^a_{bc}\, \sum_{xy=ag} \delta_{b,x}\,\delta_{c,y}\,
\xi^b\tens \xi^c\cr
&=& -\,\hat\Gamma^a_{bc}\, \xi^b\tens \xi^c\,(\delta_{g,e}-\delta_{bc,ag})\ .
\end{eqnarray*}
For this to be well defined we have to get zero when this formula is applied to
a sum of $\delta_g$ in the kernel of $\varpi$.
This gives the condition.
For $d\in C$ we have $\varpi(\delta_d)=\xi^d$, and then we have
\begin{eqnarray*}
\hat\sigma(\xi^a\tens\xi^d) &=&
 \delta_{bc,ad}\,\hat\Gamma^a_{bc}\, \xi^b\tens \xi^c\ ,
\end{eqnarray*}
and using (\ref{YDhfhdjn}) gives the formula.\quad$\square$

\begin{propos} \label{kzsbcvkcvvv}
The condition for $\nabla$
 to preserve the metric is that the matrix $g_{a,b}\Gamma^b_c$ (summation over $b$) is antiHermitian. If $g_\bullet$ is diagonal, with all enteries on the diagonal equal
(and necessarily real), then this reduces to $\hat\Gamma^a_{d,c}\,=\,(\hat\Gamma^{c}_{d^{-1},a})^*$.
\end{propos}
\proof The matrix $P$ is the identity, so $\extd P=0$, and $g_\bullet $ is complex valued, so again $\extd g_\bullet=0$. Using \ref{ksdzjcv}, this gives the equation
\begin{eqnarray*}
g_{a,b}\,\hat\Gamma^b_{d,c}\,=\,g_{c,b}\,(\hat\Gamma^{b}_{d^{-1},a})^*\ .\quad\square
\end{eqnarray*}

\begin{propos} \label{sfkhjgasjhk} For $(\nabla,\sigma)$ a bimodule covariant derivative
as in \ref{iusvggiuav},
$\nabla$ is torsion compatible if and only if, for all
$b,c,d\in C$,
\begin{eqnarray*}
d^{-1}bc\in C\quad \Longrightarrow \quad
\hat\Gamma^d_{b,c}-\hat\Gamma^d_{c,c^{-1}bc} &=& \delta_{cd,bc}-\delta_{b,d}\ .
\end{eqnarray*}
\end{propos}
\proof\quad We use (\ref{iusvggiuav}) to write
\begin{eqnarray*}
(\id+\sigma)(\xi^d\tens\xi^k) &=&  \delta_{bc,dk}\,\hat\Gamma^d_{bc}\, \xi^b\tens \xi^c
+ \xi^{dkd^{-1}}\tens \xi^d + \xi^d\tens\xi^k\ ,
\end{eqnarray*}
and asking that this is the same as $\Psi$ applied to it yields
\begin{eqnarray*}
\delta_{bc,dk}\,\hat\Gamma^d_{bc}\, \xi^b\tens \xi^c
+  \xi^d\tens\xi^k
&=&
\delta_{bc,dk}\,\hat\Gamma^d_{bc}\, \xi^{bcb^{-1}}\tens \xi^b
+ \xi^{dk({dkd^{-1}})^{-1}}\tens \xi^{dkd^{-1}}  \ .
\end{eqnarray*}
Rearranging this gives
\begin{eqnarray*}
\delta_{bc,dk}\,(\hat\Gamma^d_{b,c}-\hat\Gamma^d_{c,c^{-1}bc})\, \xi^b\tens \xi^c
&=&
 \xi^{dk({dkd^{-1}})^{-1}}\tens \xi^{dkd^{-1}} -\xi^d\tens\xi^k \ ,
\end{eqnarray*}
and this can be
rewritten as
\begin{eqnarray*}
\delta_{bc,dk}\,(\hat\Gamma^d_{b,c}-\hat\Gamma^d_{c,c^{-1}bc})
&=& \delta_{b,dk({dkd^{-1}})^{-1}}\, \delta_{c,dkd^{-1}}-\delta_{b,d}\, \delta_{c,k}\ .
\end{eqnarray*}
Observe that if we have $bc\neq dk$ we get $0=0$ for this equation.
Then we set $k=d^{-1}bc$, giving
\begin{eqnarray*}
d^{-1}bc\in C\quad \Longrightarrow \quad \hat\Gamma^d_{b,c}-\hat\Gamma^d_{c,c^{-1}bc}
&=& \delta_{b,bc({bcd^{-1}})^{-1}}\, \delta_{c,bcd^{-1}}-\delta_{b,d}\, \delta_{c,d^{-1}bc}\cr
&=&  \delta_{cd,bc}-\delta_{b,d}\ .\quad\square
\end{eqnarray*}

\begin{propos} \label{jkxhbvbvc}
The condition for star compatibility (see \ref{badsuvcvbi}) to hold is, summing over $b'$,
\begin{eqnarray*}
c^{-1}ab\in C\quad\Longrightarrow \quad
(\hat\Gamma^a_{abb'^{-1},b'}+\delta_{a,b'})
  ((\hat\Gamma^{b'^{-1}}_{b^{-1}a^{-1}c,c^{-1}})^*+\delta_{b'^{-1},c^{-1}}) \,=\,\delta_{a,c} \ .
\end{eqnarray*}
\end{propos}
\proof Using the formula for $\sigma_{\overline{\Omega^1}}$ from \ref{sliogcb},
we have
\begin{eqnarray*}
\sigma_{\overline{\Omega^1}}(\star\tens\id)
 &=&  (\star^{-1}\tens\id) \Upsilon\, \overline{\sigma_{\Omega^1}^{-1}}\, \Upsilon^{-1}\,
(\star\tens\star)\ .
\end{eqnarray*}
so the condition we want is
\begin{eqnarray*}
(\star\tens\star) \, \sigma_{\Omega^1}
 &=&  \Upsilon\, \overline{\sigma_{\Omega^1}^{-1}}\, \Upsilon^{-1}\,
(\star\tens\star)\ .
\end{eqnarray*}
In our present case,
\begin{eqnarray*}
\overline{\sigma_{\Omega^1}}\,\Upsilon^{-1}\,(\star\tens\star) \, \sigma_{\Omega^1}(\xi^a\tens\xi^b)
 &=&   \Upsilon^{-1}\,
(\star\tens\star)(\xi^a\tens\xi^b)\ ,
\end{eqnarray*}
and this gives (summing over all primed indices)
\begin{eqnarray*}
 \overline{\xi^{b^{-1}}\tens\xi^{a^{-1}}} &=& \delta_{ab,a'b'}(\hat\Gamma^a_{a',b'}+\delta_{a,b'})
\overline{\sigma_{\Omega^1}}\,\Upsilon^{-1}\,(\star\tens\star) (\xi^{a'}\tens\xi^{b'}) \cr
 &=& \delta_{ab,a'b'}(\hat\Gamma^a_{a',b'}+\delta_{a,b'})
\overline{\sigma_{\Omega^1}(\xi^{b'^{-1}} \tens \xi^{a'^{-1}})} \cr
 &=& \delta_{ab,a'b'}(\hat\Gamma^a_{a',b'}+\delta_{a,b'})\delta_{b'^{-1}a'^{-1},b''a''}
 (\hat\Gamma^{b'^{-1}*}_{b'',a''}+\delta_{b'^{-1},a''})\,  \overline{\xi^{b''}\tens\xi^{a''}}\ ,
\end{eqnarray*}
and this gives (where the sum is over $a'$ and $b'$)
\begin{eqnarray*}
\delta_{a^{-1},a''}\, \delta_{b^{-1},b''}
 &=& \delta_{ab,a'b'}\, \delta_{b^{-1}a^{-1},b''a''}\,(\hat\Gamma^a_{a',b'}+\delta_{a,b'})
  (\hat\Gamma^{b'^{-1}*}_{b'',a''}+\delta_{b'^{-1},a''})\ .
\end{eqnarray*}
This reduces to (summing over $b'$)
\begin{eqnarray*}
\delta_{a,c}\, \delta_{b,d}
 &=& \delta_{ab,cd}\,(\hat\Gamma^a_{abb'^{-1},b'}+\delta_{a,b'})
  (\hat\Gamma^{b'^{-1}*}_{d^{-1},c^{-1}}+\delta_{b'^{-1},c^{-1}}) \ ,
\end{eqnarray*}
and substituting for $d$ gives  (summing over $b'$)
\begin{eqnarray*}
\delta_{a,c}\, \delta_{b,d}
 &=& \delta_{ab,cd}\,(\hat\Gamma^a_{abb'^{-1},b'}+\delta_{a,b'})
  (\hat\Gamma^{b'^{-1}*}_{b^{-1}a^{-1}c,c^{-1}}+\delta_{b'^{-1},c^{-1}}) \ .
\end{eqnarray*}
If we choose $a\neq c$ we see that  (summing over $b'$)
\begin{eqnarray*}
c^{-1}ab\in C\quad\Longrightarrow \quad
(\hat\Gamma^a_{abb'^{-1},b'}+\delta_{a,b'})
  (\hat\Gamma^{b'^{-1}*}_{b^{-1}a^{-1}c,c^{-1}}+\delta_{b'^{-1},c^{-1}})\,=\,0 \ ,
\end{eqnarray*}
wheras $a=c$ gives
\begin{eqnarray*}
(\hat\Gamma^a_{abb'^{-1},b'}+\delta_{a,b'})
  (\hat\Gamma^{b'^{-1}*}_{b^{-1}a^{-1}c,c^{-1}}+\delta_{b'^{-1},c^{-1}})\,=\,1 \ .\quad\square
\end{eqnarray*}

\subsection{Bimodule covariant derivatives on the permutation group $S_3$}
We take the example where $G=S_3$, the permutations of three objects, and $C$ to
be the subset of transpositions.
As the product of any three elements of $C$ is in $C$, the condition for the
existence of $\sigma$ in \ref{iusvggiuav} gives no
restrictions on the numbers
$\hat\Gamma^a_{bc}\in\mathbb{C}$, so any covariant derivative is a bimodule covariant derivative.

If $G$ is a right module map, then $g^{a,b}$ has only diagonal enteries non-zero.
As $C$ is a single conjugacy class, invariance of $G$ to the right coaction corresponds to
all the diagonal enteries of $g^{a,b}$ being the same. To be concrete, we will use this metric for the rest of the example.

Similarly, to simplify matters, we will restrict attention to $\nabla$ being right invariant. By conjugation of the indices, we see that there are
only 5 possible different values of the Christoffel symbols $\hat\Gamma^a_{bc}$, which we call
$\hat\Gamma^x_{xx}$, $\hat\Gamma^x_{yz}$, $\hat\Gamma^x_{yx}$,
$\hat\Gamma^x_{xy}$ and $\hat\Gamma^y_{xx}$ (where $x,y,z$ are all different). We set
\begin{eqnarray}\label{bvdhuas}
\hat\Gamma^x_{xx}=a-1\ ,\quad \hat\Gamma^x_{yz}=c\ ,\quad \hat\Gamma^x_{yx}=d-1\ ,\quad
\hat\Gamma^x_{xy}=e\ ,\quad \hat\Gamma^y_{xx}=b\ .
\end{eqnarray}
We write $1,2,3$ for the transpositions, and use this to order the set $C$.
Now we shall build up the matrix for $\sigma$, using the convention
\begin{eqnarray*}
A\tens B \,=\, \left(\begin{array}{ccc}A_{11}\, B & A_{12}\, B & A_{13}\, B
\\ A_{21}\, B & A_{22}\, B & A_{23}\, B
\\ A_{31}\, B & A_{32}\, B & A_{33}\, B
\end{array}\right)
\end{eqnarray*}
and a little calculation gives
\begin{eqnarray*}
\sigma \,=\,\left(
\begin{array}{lllllllll}
 a & 0 & 0 & 0 & b & 0 & 0 & 0 & b \\
 0 & e & 0 & 0 & 0 & d & c & 0 & 0 \\
 0 & 0 & e & c & 0 & 0 & 0 & d & 0 \\
 0 & 0 & d & e & 0 & 0 & 0 & c & 0 \\
 b & 0 & 0 & 0 & a & 0 & 0 & 0 & b \\
 0 & c & 0 & 0 & 0 & e & d & 0 & 0 \\
 0 & d & 0 & 0 & 0 & c & e & 0 & 0 \\
 0 & 0 & c & d & 0 & 0 & 0 & e & 0 \\
 b & 0 & 0 & 0 & b & 0 & 0 & 0 & a
\end{array}
\right)
\end{eqnarray*}
and hence
\begin{equation}\label{detsigma}
\det(\sigma)=(a-b)^2 (a+2 b) (e+c+d)^2 \left(e^{2}-c e-d e+c^2+d^2-c d\right)^2.\end{equation}

\begin{lemma}\label{lembras3} The covariant derivative $\nabla$ in (\ref{bvdhuas}) has invertible  $\sigma$ obeying the braid relations if and only if one of the  following:

1)\quad $b=c=d=0$, $a=e\neq 0$

2)\quad $b=e=c=0$, $d=a\neq 0$

3) \quad $b=d=e=0$, $c=a\neq 0$

4)\quad $b\neq 0$, $c=d=0$, $e=-a^2/b=a-b$

5)\quad  $b,c\neq 0$, $c=d=-(b+e)$, $a=-eb/(b+e)$, $b^2+be+e^2=0$

6)\quad $b,c\neq 0$, $a=c=e$, $d=b^2/e$, $b^2+be+e^2=0$

7)\quad $b,c\neq 0$, $a=e$, $c=-b-e$, $d=-b^2/(b+e)$, $b^2+be+e^2=0$
\end{lemma}
\proof Mathematica calculation.\quad$\square$

\medskip
Now we consider star compatibility, in the form given in
\ref{jkxhbvbvc}. The easiest way to do this is to consider the matrices
 $\hat\Gamma^a_{b,c}+\delta_{a,c}=N(bc)_{a,c}$, which are in this case
\begin{eqnarray} \label{ksjbvbvbcv}
N(e) &=& \left(\begin{array}{ccc}a & b & b \\b & a & b \\b & b & a\end{array}\right)\ ,\cr
N(12) &=& \left(\begin{array}{ccc}d & e & c \\c & d & e \\ e & c & d\end{array}\right)\ ,\cr
N(21) &=& \left(\begin{array}{ccc}d & c & e \\ e & d & c \\c & e & d\end{array}\right)\ .
\end{eqnarray}

\begin{propos} \label{bvhsdiavbas}
The covariant derivative $\nabla$ as in (\ref{bvdhuas}) is star compatible, in the form given in
\ref{jkxhbvbvc}, if and only if the matrices $N(e)$ and $N(12)$ in (\ref{ksjbvbvbcv}) are unitary.
\end{propos}
\proof Remember that in our present case, $a=a^{-1}$ for all $a\in C$.
Then the conclusion of \ref{jkxhbvbvc} can be viewed as a matrix multiplication. Also
use the fact that $N(e)^T=N(e)$ and $N(12)^T=N(21)$. \quad$\square$

\begin{propos}
The covariant derivative $\nabla$  in (\ref{bvdhuas}) preserves the metric if and only if
$a,c,d$ are real and $e=b^*$.
\end{propos}
\proof
{} From \ref{kzsbcvkcvvv}. \quad$\square$

\begin{propos}\label{bcadhsudbvcaui}
The covariant derivative $\nabla$  in (\ref{bvdhuas}) is torsion compatible
(see \ref{sfkhjgasjhk}) if and only if  it is torsion free, and this holds if and only if $d=c=e+1$.
\end{propos}
\proof We first look at torsion compatibility.
From \ref{sfkhjgasjhk} we have
\begin{eqnarray*}
\hat\Gamma^d_{b,c}-\hat\Gamma^d_{c,c^{-1}bc} &=& \delta_{cd,bc}-\delta_{b,d}\ ,
\end{eqnarray*}
and we consider the cases.

Assuming that we have torsion compatibility, and using \ref{sfkhjgasjhk}, we have
\begin{eqnarray*}
-\wedge\nabla\xi^x &=& \sum_{y\neq x}\Big( \hat\Gamma^x_{xy}\,\xi^x\wedge\xi^y
+ \hat\Gamma^x_{yx}\,\xi^y\wedge\xi^x\Big)
+ \sum_{y,z:\,x,y,z\ \mathrm{different}} \hat\Gamma^x_{yz}\,\xi^y\wedge\xi^z \cr
 &=& \sum_{y\neq x}\Big( e\,\xi^x\wedge\xi^y
+(d-1)\,\xi^y\wedge\xi^x\Big)
+c\, \sum_{y,z:\,x,y,z\ \mathrm{different}} \xi^y\wedge\xi^z \ ,
\end{eqnarray*}
hence
\begin{eqnarray*}
\extd\xi^x-\wedge\nabla\xi^x  &=& \sum_{y\neq x}\Big( (e+1)\,\xi^x\wedge\xi^y
+ d\,\xi^y\wedge\xi^x\Big)
+c\, \sum_{y,z:\,x,y,z\ \mathrm{different}} \xi^y\wedge\xi^z \ .
\end{eqnarray*}
Setting $x=1$ say, and using the relations in $\Omega^2$, which are
\begin{eqnarray*}
&&\xi^1\wedge\xi^1=\xi^2\wedge\xi^2=\xi^3\wedge\xi^3=0\ ,\cr
&& \xi^1\wedge\xi^2+\xi^2\wedge\xi^3+\xi^3\wedge\xi^1=0\ ,\quad
 \xi^2\wedge\xi^1+\xi^3\wedge\xi^2+\xi^1\wedge\xi^3=0\ ,
\end{eqnarray*}
we expand in a basis $\xi^1\wedge\xi^2$, $\xi^1\wedge\xi^1$, $\xi^2\wedge\xi^3$ and
$\xi^3\wedge\xi^2$ say, to conclude that the torsion vanishes if
$d=c=e+1$. Hence torsion compatible implies torsion free. The converse is immediate.
\quad$\square$

\begin{propos}\label{cvcvictcttcttcgg}
The covariant derivative $\nabla$  in (\ref{bvdhuas}) has
 vanishing cotorsion (see definition \ref{kxzcvzvck}) if and only if,
\begin{eqnarray*}
 e-b^*\,=\, c-c^*
\,=\, d-d^*\
\end{eqnarray*}
which is weaker than full metric preservation.
\end{propos}
\proof
  From \ref{cvcvictcttcttc}, using the fact that $\extd P=0$ and $\extd g_\bullet=0$, we get
\begin{eqnarray*}
(g_\bullet \Gamma+\Gamma^* g_\bullet)
\wedge e^\bullet &=& 0\ .
\end{eqnarray*}
As $g_\bullet$ is a multiple of the identity matrix this gives
$(\Gamma+\Gamma^*)\wedge e^\bullet=0$, or inserting the indices again,
\begin{eqnarray*}
(\Gamma^i_j+(\Gamma^j_i)^*)\wedge \xi^j\,=\,0
\end{eqnarray*}
which is
\begin{eqnarray*}
(\hat\Gamma^i_{kj}\,\xi^k+(\Gamma^j_{ki}\,\xi^k)^*)\wedge \xi^j\,=\,0\ ,
\end{eqnarray*}
which becomes the following equation,  for all $i$:
\begin{eqnarray*}
(\hat\Gamma^i_{kj}\,\xi^k-(\hat\Gamma^j_{ki})^*\,\xi^{k^{-1}})\wedge \xi^j\,=\,0\ ,
\end{eqnarray*}
(in our case $k^{-1}=k$).
Now put $i=1$, and using the fact that $\xi^j\wedge\xi^j=0$, we get
\begin{eqnarray*}
0 &=&   (\hat\Gamma^1_{12}-(\hat\Gamma^2_{11})^*)\,\xi^1\wedge \xi^2 +
(\hat\Gamma^1_{13}-(\hat\Gamma^3_{11})^*)\,\xi^1\wedge \xi^3   \cr
&& +\, (\hat\Gamma^1_{21}-(\hat\Gamma^1_{21})^*)\,\xi^2\wedge \xi^1
+  (\hat\Gamma^1_{23}-(\hat\Gamma^3_{21})^*)\,\xi^2\wedge \xi^3    \cr
&& +\, (\hat\Gamma^1_{31}-(\hat\Gamma^1_{31})^*)\,\xi^3\wedge \xi^1
+  (\hat\Gamma^1_{32}-(\hat\Gamma^2_{31})^*)\,\xi^3\wedge \xi^2\ ,
\end{eqnarray*}
and using the explicit relations in $\Omega^2$ as in the preceeding proof, the conditions for this are
\begin{eqnarray*}
&& \hat\Gamma^1_{12}-(\hat\Gamma^2_{11})^*\,=\, \hat\Gamma^1_{23}-(\hat\Gamma^3_{21})^*
\,=\, \hat\Gamma^1_{31}-(\hat\Gamma^1_{31})^*\ ,\cr
&& \hat\Gamma^1_{13}-(\hat\Gamma^3_{11})^* \,=\, \hat\Gamma^1_{21}-(\hat\Gamma^1_{21})^*
\,=\, \hat\Gamma^1_{32}-(\hat\Gamma^2_{31})^*\ .
\end{eqnarray*}
Substituting the values $a,b,c,d,e$ here gives
\begin{eqnarray*}
&& e-b^*\,=\, c-c^*
\,=\, d-d^*\ ,\cr
&& e-b^* \,=\, d-d^*
\,=\, c-c^*\ . \quad\square
\end{eqnarray*}

 \begin{propos} \label{jdhvzbks}
The covariant derivative $\nabla$  in (\ref{bvdhuas}) is torsion free and is star compatible if and only if $d=c$, $e=c-1$  and
\[  a={e^{\imath\phi}\over \sqrt{1+8\,\cos^2(\theta-\phi)}},\quad b=-{e^{\imath\theta}\, 2\,\cos(\theta-\phi)\over  \sqrt{1+8\,\cos^2(\theta-\phi)}},\quad c={2 e^{\imath\psi}\cos\psi\over 3}\]
for arbitrary angles subject to  $|\theta-\phi|, |\psi|\le {\pi\over 2}$ and certain identifications on the boundary (a 3-parameter moduli space). \end{propos}
\proof
The unitarity of $N(e)$ (see \ref{bvhsdiavbas}) requires
\begin{eqnarray*}
|a|^2+2\,|b|^2\,=\,1\ ,\quad |b|^2+a^*\,b+a\,b^*\,=\,0
\end{eqnarray*}
which we solve for  $a=R\,e^{\mathrm{i}\,\phi}$
and $b=-r\,e^{\mathrm{i}\,\theta}$ where $r,R\ge 0$ as stated (the - sign here is for convenience). The first equation tells us that $R=\sqrt{1-2r^2}$ while the second for $r>0$ tells us that $r=\sqrt{1-2r^2}2\cos(\theta-\phi)$ which provides us $r=2\cos(\theta-\phi)/\sqrt{1+8\cos^2(\theta-\phi)}$ and positive cosine, i.e. $|\theta-\phi|<{\pi\over 2}$. If $r=0$ then the second equation is empty, $R=1$ and $\phi$ is unconstrained. We think of this as on the boundary $\theta=\phi\pm {\pi\over 2}$ with the two values of $\theta$ identified.  The conditions for torsion compatibility imply that $N(12)=c\,A-B$, where
$A$ and $B$ are the matrices
\begin{eqnarray*}
A\,=\, \left(\begin{array}{ccc}1 & 1 & 1 \\1 & 1 & 1 \\1 & 1 & 1\end{array}\right)\ ,\quad
B\,=\, \left(\begin{array}{ccc}0 & 0 & 1 \\1 & 0 & 0 \\0 & 1 & 0\end{array}\right)\ .
\end{eqnarray*}
Then we see that $N(12)^*\, N(12)=(3\,|c|^2 -c-c^*)A+I_3$. This unitarity condition therefore amounts
to
\begin{equation}\label{ceqn}   3\,|c|^2\,=\,c+c^*\end{equation}
which we solve for $c=\rho e^{\imath\psi}$ as $\rho=0$ or $3\rho=2\cos\psi$ and $|\psi|<{\pi\over 2}$. We write the first as the identified boundaries of the open interval for $\psi$. This describes our parametrization of the moduli space of such connections. \quad$\square$

\begin{propos} \label{bchadsbcv}
The moduli space of covariant derivatives $\nabla$  in (\ref{bvdhuas}) which are torsion free, cotorsion free and star compatible has one continuous parameter $r\in[{1\over 3},{2\over 3}]$ and connection as in Proposition~\ref{jdhvzbks} with
\[ \cos\theta={1+3r^2\over 4r},\quad \cos(\theta-\phi)={r\over 2\sqrt{1-2r^2}},\quad \cos\psi=\sqrt{{9\over 8}(1-r^2)}\]
and free choices for the sign of $\theta,\psi,\theta-\phi$.
\end{propos}
\proof From \ref{cvcvictcttcttcgg}
the condition for vanishing cotorsion is that
\begin{eqnarray*}
 e-b^*\,=\, c-c^*
\,=\, d-d^*\ .
\end{eqnarray*}
Combining this with the vanishing torsion condition $c=d=e+1$ gives $c=1+b$ as the content of the cotorsion free condition if the torsion is known to vanish. Hence we merely need to constrain  the moduli space  in Proposition~\ref{jdhvzbks} by this requirement:
\[ -r\sin\theta={2\over 3}\cos(\psi)\sin(\psi),\quad 1-r\cos\theta={2\over 3}\cos^2\psi\]
where $r=2\cos(\theta-\phi)/\sqrt{1+8\cos^2(\theta-\phi)}$ is the value of $|b|$. Inverting this relationship (or going back to the derivation) gives the middle displayed equation. Also by squaring the second equation it readily follows that
\begin{eqnarray*}
3\,r^2-4\,r\,\cos\theta+1\,=\,0\
\end{eqnarray*}
(this also follows at once from $c=b+1$ in (\ref{ceqn})),  providing the first displayed equation stated. Finally, rearranging, squaring and adding gives
\[ r^2=1-{8\over 9}\cos^2\psi\]
which provides the last displayed equation stated. To have solutions for $\psi$ we need $r\in[{1\over 3},1]$ and for $\theta-\phi$ we need $r\in[0,{2\over 3}]$. Given the first restriction, the equation for $\theta$ does not constrain $r$. At the endpoint $r=1/3$ we have $\psi=0$ and at $r=2/3$ we have $\theta-\phi=0$ so no choice of their signs at the respective endpoints. \quad$\square$

\medskip According to the terminology introduced in \cite{Ma:rief} such connections could be called `generalised Levi-Civita' in the sense that only vanishing cotorsion, which is weaker than metric compatibility, is required. Whereas \cite{Ma:rief} introduced a theory of frame bundles and spin connections and in that context (which did not consider $*$-structures) there was a unique `generalised Levi-Civita' connection on $S_3$ for the Euclidean metric, we see that our theory of bimodule linear connections is less restrictive, even after we introduce $*$-compatibility.  In retrospect, this should not have been completely unexpected. Complex analytic functions have isolated zeros, but in the case here the introduction of complex conjugates allows the possibility of non isolated solutions for the complex parameters in the covariant derivative.

The reader may now ask whether this multiplicity may be reduced to a unique solution if we strengthen the condition of cotorsion free to preserving the metric. The answer is that we now get no solutions at all that are also torsion free and star compatible:

\begin{propos}\label{metricstars3}
The covariant derivative $\nabla$  in (\ref{bvdhuas}):  \smallskip

\noindent\textbf{1)} preserves the metric and is torsion free if and only if
$a,c$ are real, $c=d=e+1=b+1$, i.e. a 2-parameter moduli space of such connectons.

\smallskip\noindent\textbf{2)} preserves the metric and is star compatible if and only if it is one of the following discrete moduli of possibilities:

i)\quad $b=e=c=0$, $a=\pm1$, $d=\pm 1$ (independent signs)

ii)\quad $b=e=d=0$, $a=\pm1$, $c=\pm 1$ (independent signs)

iii)\quad $a=d=\pm\frac13$, $b=e=c=\mp\frac23$

iv)\quad $a=c=\pm\frac13$, $b=e=d=\mp\frac23$.

\smallskip\noindent\textbf{3)} cannot be metric preserving, torsion free and star compatible.

\end{propos}
\proof
We recall that metric preserving requires  that $a,c,d$ are real and $e=b^*$.
Part (1) is then immediate. For (2) we have the unitarity requirements:
\begin{eqnarray*}
\left(\begin{array}{ccc}a & b & b \\b & a & b \\b & b & a\end{array}\right)\
\left(\begin{array}{ccc}a & b^* & b^* \\b^* & a & b^* \\b^* & b^* & a\end{array}\right) &=&
\left(\begin{array}{ccc}1 & 0 & 0 \\0 & 1 & 0 \\0 & 0 & 1\end{array}\right)\ ,\cr
 \left(\begin{array}{ccc}d & b^* & c \\c & d & b^* \\b^* & c & d\end{array}\right)\
  \left(\begin{array}{ccc}d & c & b \\b & d & c \\c & b & d\end{array}\right)&=&
\left(\begin{array}{ccc}1 & 0 & 0 \\0 & 1 & 0 \\0 & 0 & 1\end{array}\right)\ .
\end{eqnarray*}
This gives the following equations:
\begin{eqnarray*}
a^2+2\,|b|^2\,=\,1\ ,\quad a(b+b^*)+|b|^2\,=\,0\ ,\quad d^2+c^2+|b|^2\,=\,1\ ,\quad
dc+b^*d+cb\,=\,0\ .
\end{eqnarray*}
If $b=0$ we immediately get the first two cases. For part (3), it is easy to see that parts
(1) and (2) have no intersection.

The assumption that $b\neq0$ gives  the following from the first two equations
for $b=x+\mathrm{i}\,y$: Both $x$ and $a$ are nonzero and
\begin{eqnarray*}
x\,=\,\frac{a^2-1}{4\,a}\ ,\quad y^2\,=\, -\,\frac{(a-1)(a+1)(3a-1)(3a+1)}{16\,a^2}\ .
\end{eqnarray*}
If we assume $y\neq 0$, then we find $d=c=-2x$, and then $9\,x^2+y^2=1$, which is not consistent with the above. We deduce that $y=0$, so $a=\pm\frac13$. Then $b=\mp\frac23$
and we are left with
\begin{eqnarray*}
c^2+d^2\,=\,\frac59\ ,\quad c\,d\,=\, \pm\frac23\, (c+d)\ .\quad\square
\end{eqnarray*}

Ironically, the  unique torsion and cotorsion free connection on $C(S_3)$ found in the frame bundle approach of \cite{Ma:rief} has covariant derivative is given by $\nabla$  in (\ref{bvdhuas})
with parameters
$a=\frac53$, $d=c=\frac23$ and $e=b=-\frac13$.
In terms of the Hermitian metric theory in the present paper (incuding the bar), this means that it is part of our 2-parameter family that is torsion free and actually preserves the Hermitian
metric, but it is not star compatible.

A priori, any one of the discrete moduli in the above part 2) (i.e., metric preserving and star compatible connections) could be considered the strictly metric preserving  `Levi-Civita' one for the Euclidean metric that we have adopted. They all have invertible $\sigma$ and we see that they all have torsion, which therefore seems to be forced in the finite theory  if we require the metric to be preserved exactly.
Finally, we see that imposing the braid relations can pin these down further.

\begin{propos} The covariant derivative $\nabla$ in  (\ref{bvdhuas}) has $\sigma$ which:

\smallskip\noindent{\bf 1)} is invertible and obeys the braid relations, and $\nabla$ is torsion free  if and only if the connection belongs to the discrete moduli of possibilities

i) $a=e=-1$, $b=c=d=0$

ii) $c=d=0$, $e=-1$, $b={1\over 2}(3\pm\sqrt{5})$, $a=b-1$

iii) $a={ 3\pm\imath{\sqrt{3}}\over 3\pm\imath 3\sqrt{3}}$, $b=\pm{ \imath\over \sqrt{3}}$, $e=-{-3\pm \imath\sqrt{3}\over 6}$, $c=d={3\mp\imath\sqrt{3}\over 6}$

\smallskip\noindent{\bf 2)} is invertible and obeys the braid relations, and $\nabla$ preserves the metric if and only if the connection belongs to one of the 1-parameter moduli components

i) $b=e=c=0$, $d=a\ne 0$ is real

ii) $b=d=e=0$, $c=a\ne 0$ is real

iii)  $c=d=a\ne0$ is real, $e=ae^{\pm{\imath2\pi\over 3}}=b^*$.

\smallskip\noindent{\bf 3)} cannot be invertible and obey the braid relations with $\nabla$  torsion free and either  cotorsion free or star compatible.
\end{propos}
\proof We conjunct our earlier results. For part 1) of the statement, of the possibilities 1)-7) for the braiding in Lemma~\ref{lembras3} only 1),4),5) have joint solutions with $c=d=e+1$, giving i),ii),iii) respectively. None of these have $b=e$ so cotorsion free is excluded in this context. Also, none of them have the matrices (\ref{ksjbvbvbcv}) unitary, giving part 3) of the statement. For part 2) of the statement, of the possibilities for the braiding only cases 2),3),5) in the classification have joint solutions with $a,c,d$ real and $b=e^*$, giving the cases i),ii) and iii) in this part of the statement respectively. The latter case requires us to solve $e^*{}^2+e^*e+e^2=0$, which then determines the rest as $b=e^*$,  $c=d=-(e+e^*)$ and $a=-ee^*/(e+e^*)=c$. We solve for $e$ in polar form. \quad$\square$

\begin{corol}\label{corbras3}  The covariant derivative $\nabla$ in (\ref{bvdhuas}) is metric preserving, star compatible and has invertible $\sigma$ obeying the braid relations if and only if it is one of the following discrete moduli of possibilities:

i)\quad  $b=e=c=0$, $a=d=\pm1$

ii)\quad $b=d=e=0$, $c=a=\pm 1$.

\noindent All cases necessarily have torsion and $\sigma^3=\pm1$ according to the sign. The +1 choice in i) is the left Maurer-Cartan connection $\nabla^L=0$.
\end{corol}
\proof Comparing the last result with Proposition~\ref{metricstars3}, we see that the only conjunction with metric and star compatibility are cases 2)i) or 2)ii) in the last proposition.   The $+1$ case corresponds to $\hat\Gamma=0$ and hence to $\nabla^L=0$. \quad$\square$

Thus, the standard Maurer-Caratan connection can be viewed as  `Levi-Civita--with-- torsion', in the sense of metric preserving and $*$-compatible, and there is also a further `non-standard' choice (it is not the right Maurer-Cartan form as this has the inverse braiding and hence $\sigma^3=1$ again).
\bigskip

\section{Example: Riemannian geometry on quantum $SU(2)$ with the 3D calculus}
Suppose that $q\in\mathbb{C}$ with $q^2\ne 1$.
The quantum group $\C_q[SL_2]$ has generators $a,b,c,d$
 with relations:
\[ ba=qab\ ,\ ca=qac\ ,\ db=qbd\ ,\ dc=qcd\ , \ cb=bc \ , \
 da-ad=q(1-q^{-2}) bc \ , \ ad-q^{-1}bc=1\]
The coproduct $\Delta$ and counit $\eps$ have the usual matrix
coalgebra form. We denote the antipode or `matrix inverse' by $S$:
\begin{eqnarray*}
S\left(\begin{array}{cc}a & b \\c & d\end{array}\right) \,=\,
\left(\begin{array}{cc}d & -q\,b \\ -q^{-1}\,c & a\end{array}\right)\ .
\end{eqnarray*}
 On $\C_q[SL_2]$ we take the 3D calculus of
\cite{Wor:dif}. In our conventions this has a basis
\[  e^-=d.\extd b-q b.\extd
d,\quad  e^+=q^{-1} a.\extd c-q^{-2}c.\extd a,\quad  e^0=d.\extd a-q
b.\extd c\]
of left-invariant 1-forms, is spanned by these as a left
module (according to the above) while the right module relations
and exterior derivative are given in these terms by:
\begin{eqnarray*}
&& e^\pm
\left(\begin{array}{cc}a & b \\c & d\end{array}\right)
=\left(\begin{array}{cc}qa & q^{-1}b \\qc & q^{-1}d\end{array}\right)
e^\pm,\quad  e^0
\left(\begin{array}{cc}a & b \\c & d\end{array}\right)
=\left(\begin{array}{cc}q^2 a & q^{-2}b \\q^2c & q^{-2}d\end{array}\right)
 e^0 \ ,\cr \cr
&& \extd a=a e^0+q b e^+,\quad \extd
b=a e^--q^{-2}b e^0,\quad \extd c=c  e^0+q d e^+,\quad \extd d=c
e^--q^{-2}d e^0
\end{eqnarray*}
For $\C_q[SL_2]$ the natural extension compatible with
the super-Leibniz rule on higher forms and $\extd^2=0$ is:
\[ \extd  e^0=q^3 e^+\wedge e^-,\quad
\extd e^\pm=\mp q^{\pm 2}[2;q^{- 2}]e^\pm\wedge e^0,\quad
e^\pm\wedge e^\pm= e^0\wedge e^0=0\]
\[
q^2 e^+\wedge e^-+ e^-\wedge e^+=0,\quad  e^0\wedge e ^\pm+q^{\pm
4}e^\pm\wedge e^0=0\] where $[n;q]=(1-q^n)/(1-q)$ denotes a
$q$-integer. This means that there are the same dimensions as
classically, including a unique top form $ e^-\wedge e^+\wedge
e^0$.

The algebra $\C_q[SL_2]$ equipped with the star operation $a^*=d$, $d^*=a$, $c^*=-q\,b$
and $b^*=-q^{-1}c$, where $q$ is real, is denoted $\C_q[SU_2]$.
Moreover the ideal corresponding to the above calculus is stable under $*$, so
$\Omega^1 \C_q[SU_2]$ is a star object. Then one has
\begin{eqnarray*}
e^{0*} &=& \extd a^*.d^*-q\extd c^*.b^* = \extd d.a-q\extd b.c =\, -e^0\ , \cr
e^{+*} &=& q^{-1}\extd c^*.a^*-q^{-2}\extd a^*.c^* = -\extd b.d+q^{-1}\extd d.b =\, -q^{-1}\, e^-\ , \cr
e^{-*} &=& \extd b^*.d^*-q\extd d^*.b^*=  \,-q\,e^+
\end{eqnarray*}
using the relations above.

The right coadjoint action of $\C_q[SU_2]$ on the left invariant one forms is particularly simple.
There is a Hopf* algebra map $\pi:\C_q[SU_2]\to \mathbb{CZ}$ (the group
algebra of the group $(\mathbb{Z},+)$, with group
generator $z$ and $z^*=z^{-1}$) given by $\pi(a)=z$,  $\pi(b)=\pi(c)=0$ and $\pi(d)=z^{-1}$.
There is a right action of $ \mathbb{CZ}$ on $\Lambda^1\C_q[SU_2]$ given by
$e^\pm\hat\ra z=q\,e^\pm$ and $e^0\hat\ra z=q^2\,e^0$,
and the right coadjoint action of $\C_q[SU_2]$ on $\Lambda^1 \C_q[SU_2]$ (written $\ra$) can be
written in terms
of $\pi$ and $\hat\ra$ as
$\xi\ra h=\xi\hat\ra\pi(h)$.

Although the calculus is not bicovariant under $\C_q[SU_2]$,
the projected right coaction $(\id\tens\pi)
\Delta:\C_q[SU_2]\to \C_q[SU_2]\tens \mathbb{CZ}$ is differentiable in the sense that we have
compatible with $\extd$ a right coaction
$\rho:\Omega^1 \C_q[SU_2]\to \Omega^1 \C_q[SU_2]\tens \mathbb{CZ}$,
which we compute as
\begin{eqnarray*}
\rho(e^0)=e^0 \tens 1\ ,\quad \rho(e^\pm)= e^\pm \tens z^{\pm 2}\ .
\end{eqnarray*}
Using $\hat\ra$ and $\rho$, $\Lambda^1$ is a right-right crossed module over $\mathbb{CZ}$ (i.e.\ a Drinfeld-Radford-Yetter module), and the corresponding braiding can be written
as (where $\pm'$ is an independent copy of $\pm$)
\begin{eqnarray}\label{cbicndiuwbc}
\Psi(e^0\tens e^0) &=& e^0\tens e^0\ , \nonumber\\
\Psi(e^0\tens e^\pm)&=& q^{\pm 4}\, e^\pm \tens e^0\ ,\nonumber\\
\Psi(e^\pm\tens e^0)&=& e^\pm \tens e^0\, \nonumber\\
\Psi(e^\pm{}'\tens e^\pm{})&=& q^{\pm 2}\,e^\pm{}\tens e^\pm{}'\ .
\end{eqnarray}

\subsection{Bimodule covariant derivatives on $\C_q[SU_2]$}
Let $\nabla$ be a covariant derivative for the calculus above.
We recall that if it is a bimodule covariant derivative, then $\sigma$ necessarily takes the form
in formula 12.
A little calculation gives
\begin{eqnarray} \label{jzxdhdhvbc77}
\extd a_{(2)}.S^{-1}(a_{(1)}) \,=\, q^{-2}\,e^0 &,&
\extd b_{(2)}.S^{-1}(b_{(1)}) \,=\, q^{-1}\, e^-\ ,\cr
\extd c_{(2)}.S^{-1}(c_{(1)}) \,=\, q^2\,e^+ &,&
\extd d_{(2)}.S^{-1}(d_{(1)}) \,=\, -e^0 \ .
\end{eqnarray}
From this, looking only at the term $\hat\sigma$ to be added in equation 14 to the
 braiding $\Psi$
 (which coincides with (\ref{cbicndiuwbc}), and which turns out to be still well-defined in this case),
 we find that $\hat\sigma$ necessarily takes the form
\begin{eqnarray} \label{ksdhvgh}
\hat\sigma(\xi\tens q^{-2}\,e^0) &=& -\nabla^L(\xi)+
(\nabla^L(\xi\hat\ra z))\hat\ra z^{-1}\ ,\cr
\hat\sigma(\xi\tens q^{-1}\, e^-) &=& 0\ ,\cr
\hat\sigma(\xi\tens q^2\,e^+) &=& 0\ ,\cr
\hat\sigma(\xi\tens (-e^0)) &=& -\nabla^L(\xi)+
(\nabla^L(\xi\hat\ra z^{-1}))\hat\ra z\ .
\end{eqnarray}
For this to be well defined we require the first and last equations to be consistent, namely:
\begin{eqnarray}\label{compatyyyy}
-q^2\,\nabla^L(\xi)+
q^2\,(\nabla^L(\xi\hat\ra z))\hat\ra z^{-1} \,=\,
\nabla^L(\xi)-
(\nabla^L(\xi\hat\ra z^{-1}))\hat\ra z\ .
\end{eqnarray}
Hence we arrive at the following proposition.

\begin{propos}\label{fghjk}
For left invariant covariant derivative $\nabla$ is a bimodule covariant derivative
if and only if,

$\nabla^L e^\pm$ is a linear combination of
$e^0\tens e^{+}$, $e^0\tens e^{-}$, $e^{+}\tens e^0$ and $e^{-}\tens e^0$,

$\nabla^L e^0$ is a linear combination of
$e^0\tens e^0$ and $e^{\pm'}\tens e^\pm$,

\noindent
In this case $\sigma$ is
\begin{eqnarray*}
\sigma(\xi\tens e^0) &=& e^0 \tens\xi
+\nabla^L(\xi)-
(\nabla^L(\xi\hat\ra z^{-1}))\hat\ra z\ ,\cr
\sigma(\xi\tens e^\pm) &=& e^\pm \tens\xi\hat\ra z^{\pm2}\ .
\end{eqnarray*}
\end{propos}
\proof Apply (\ref{compatyyyy})
and (\ref{jzxdhdhvbc}) to the different possibilities for
$\nabla^L$. \quad$\square$

\medskip If we write some of the coefficients of $\nabla^L$ in \ref{fghjk} as
\begin{eqnarray}\label{ksdvhnvdssbg}
\nabla^L(e^0) \,=\, r\, e^0\tens e^0+\alpha\ ,
\end{eqnarray}
where $r\in\mathbb{C}$ and
$\alpha$  is a linear combination of
$e^{\pm'}\tens e^\pm$, then
\begin{eqnarray}\label{vjhzcsdvchjv}
\sigma(e^0\tens e^0) &=& (1+r(1-q^2))e^0 \tens e^0\ ,\cr
\sigma(e^\pm\tens e^0) &=& e^0 \tens e^\pm
+(1-q^2)\nabla^L(e^\pm)\ .
\end{eqnarray}
 This concludes our general examination of which left
invariant connections have a well defined map $\sigma$, and the form
of that map. For the rest of this subsection we will
always assume the form of $\nabla^L$ given in \ref{fghjk} and
(\ref{ksdvhnvdssbg}).

\begin{propos} \label{kdyrnn}
The  covariant derivative
$\nabla$ in \ref{fghjk} is invariant under the right $\mathbb{CZ}$ coaction if and only if
$\alpha$ in (\ref{ksdvhnvdssbg}) and  $\nabla^L(e^\pm)$ take the form
\begin{eqnarray*}
\alpha &=& \nu\, e^+\tens e^- + \mu\,e^-\tens e^+\ ,\cr
\nabla^L(e^\pm) &=& n_\pm\, e^0\tens e^\pm + m_\pm\, e^\pm\tens e^0\ .
\end{eqnarray*}
for some  numbers $\nu,\mu,n_\pm,m_\pm$.
In this case the matrix of Christoffel symbols (using basis order $+,0,-$) is
\begin{eqnarray*}
\Gamma \,=\, -\,\left(\begin{array}{ccc}n_+\, e^0 & m_+\, e^+ & 0 \\\mu\, e^- & r\,e^0 & \nu\, e^+ \\0 & m_-\,e^- & n_-\, e^0\end{array}\right)
\end{eqnarray*}
and
\begin{eqnarray*}
\sigma(e^\pm\tens e^0) &=& (1+(1-q^2)n_\pm)\,e^0\tens e^\pm
+ (1-q^2)m_\pm\, e^\pm\tens e^0\ .
\end{eqnarray*}
\end{propos}
\proof Try all the possibilities.\quad$\square$

\medskip We will restrict our attention to bimodule covariant derivatives of this
left invariant and right $\mathbb{CZ}$ invariant
 form.
 With our previous basis order and convention for tensor product, $\sigma$ is the matrix
\begin{eqnarray*}
\sigma\,=\,\left(
\begin{array}{ccccccccc}
 q^2 & 0 & 0 & 0 & 0 & 0 & 0 & 0 & 0 \\
 0 & m_+ \left(1-q^2\right) & 0 & q^4 & 0 & 0 & 0 & 0 & 0 \\
 0 & 0 & 0 & 0 & 0 & 0 & q^2 & 0 & 0 \\
 0 & 1+n_+ \left(1-q^2\right) & 0 & 0 & 0 & 0 & 0 & 0 & 0 \\
 0 & 0 & 0 & 0 & 1+r\left(1-q^2\right) & 0 & 0 & 0 & 0 \\
 0 & 0 & 0 & 0 & 0 & 0 & 0 &1+ n_- \left(1-q^2\right)& 0 \\
 0 & 0 & \frac{1}{q^2} & 0 & 0 & 0 & 0 & 0 & 0 \\
 0 & 0 & 0 & 0 & 0 & \frac{1}{q^4} & 0 & m_- \left(1-q^2\right) & 0
   \\
 0 & 0 & 0 & 0 & 0 & 0 & 0 & 0 & \frac{1}{q^2}
\end{array}
\right)
\end{eqnarray*}

\begin{propos}\label{dvvbjj}
 The bimodule covariant derivative
$\nabla$ in \ref{kdyrnn} has $\sigma$ obeying the braid relations if and only if
 either both $m_\pm=0$ or we have two exceptional cases:
 \begin{eqnarray*}
\mathrm{Case 1:}&& n_+(1-q^2)=-1-q^{-2}\rho(r),\quad n_-=-(1+q^2)\ ,\quad m_+=q^2n_+,\quad m_-=0\ ,\cr
   \mathrm{Case 2:}&& n_+=q^{-4}(1+q^2)
   \ ,\quad n_-(1-q^2)=-1-q^2\rho(r) ,\quad m_+=0\ ,\quad m_-=q^{-2}n_-   \ ,
\end{eqnarray*}
where $\rho(r)=1+r(1-q^2)$. For both $m_\pm=0$,  $\sigma$ has eigenvectors
\begin{eqnarray*}
q^2\,e_+\tens e_-+e_-\tens e_+   &\mathrm{eigenvalue}&  1\ ,\cr
q^2\,e_+\tens e_--e_-\tens e_+   &\mathrm{eigenvalue}&  -1\ ,\cr
e_0\tens e_0   &\mathrm{eigenvalue}& \rho(r) ,\cr
q^2\,e_+\tens e_0+e_0\tens e_+\, q_+   &\mathrm{eigenvalue}&
q^2\,q_+\ ,\cr
q^2\,e_+\tens e_0-e_0\tens e_+\,q_+   &\mathrm{eigenvalue}&
-\,q^2\,q_+\ ,\cr
q^2\,q_-\, e_0\tens e_-+e_-\tens e_0  &\mathrm{eigenvalue}&
q^{-2}\,q_-\ ,\cr
q^2\,q_-\, e_0\tens e_--e_-\tens e_0   &\mathrm{eigenvalue}&
-\,q^{-2}\,q_-\ ,\cr
e_+\tens e_+  &\mathrm{eigenvalue}&  q^2\ ,\cr
e_-\tens e_-  &\mathrm{eigenvalue}&  q^{-2}\ .
\end{eqnarray*}
where $q_\pm=\sqrt{1+n_\pm(1-q^2)}$.
\end{propos}
\proof Direct calculation with the matrix $\sigma$ using
Mathematica. \quad$\square$

\begin{propos} \label{bcidasdbvcla}
 The bimodule covariant derivative
$\nabla$ in \ref{kdyrnn} is torsion compatible if and only if
 \begin{eqnarray*}
m_+\,=\,q^4\,n_+-(1+q^2)\ ,\quad m_-\,=\, q^{-4}(n_-+1+q^2)\ .
\end{eqnarray*}
\end{propos}
\proof Mathematica calculation.\quad$\square$

\begin{propos} \label{cgvvcsuvc}
The torsion for the covariant derivative in \ref{kdyrnn} is given by
\begin{eqnarray*}
T(e^\pm) &=&  (m_\pm - q^{\pm 4}\,n_\pm  \pm q^{\pm2}(1+q^{-2})  )\, e^\pm\wedge e^0
\ ,\cr
T(e^0) &=& (\nu-q^2\,\mu-q^3)\,e^+\wedge e^-\ .
\end{eqnarray*}
Thus if the covariant derivative is torsion compatible, then $T(e^\pm)=0$.
\end{propos}
\proof Begin with
\begin{eqnarray*}
\wedge\nabla\xi^\pm &=& n_\pm\, e^0\wedge e^\pm + m_\pm\, e^\pm\wedge e^0\cr
&=& (m_\pm - q^{\pm 4}\,n_\pm    )\, e^\pm\wedge e^0\ ,\cr
\wedge\nabla\xi^0 &=& r\,e^0\wedge e^0+\nu\,e^+\wedge e^-+\mu\,e^-\wedge e^+\cr
 &=& \nu\,e^+\wedge e^-+\mu\,e^-\wedge e^+\cr
  &=& (\nu-q^2\,\mu)\,e^+\wedge e^-\ .
\end{eqnarray*}
Then use $T(e^i) = \wedge\nabla e^i-\extd e^i$. \quad$\square$

\begin{propos} \label{bchdsucvasu}
The covariant derivative in \ref{kdyrnn} is star-compatible if and only if
\begin{eqnarray*}
n_-^*\,=\,\frac{-\, n_+}{1+n_+(1-q^2)} &,& m_-^*\,=\,\frac{-\, q^{-4}\,m_+}{1+n_+(1-q^2)},\quad
r^*+r+|r|^2(1-q^2)\,=\,0\ .
\end{eqnarray*}
and is star-preserving in the sense of (\ref{bcihdbbdiuv}) if and only if $\nu\,=\,-\,\mu^*\,q^2$ also holds.
\end{propos}
\noindent {\bf Proof:}\quad Mathematica calculation.\quad$\square$

\subsection{Metric compatibility on $\C_q[SU_2]$}
The condition that a symmetric Hermitian metric $g$ is $\mathbb{CZ}$ right invariant
is that $g^{ij}=0$ for $i\neq j$.
To see this, we use $G(\overline{e^i})=e_j\,g^{ji}$, and the conditions
$e_\pm\hat\ra z=q^{-1}\,e^\pm$, $e_0\hat\ra z=q^{-2}\,e_0$, $\rho(e_0)=e_0\tens 1$, and $\rho(e_\pm)=e_\pm\tens z^{\mp 2}$. In line with our restriction on $\nabla$, we will restrict attention to such $\mathbb{CZ}$ right invariant
metrics. We write $g$ as a matrix $g_\bullet$ in the basis $e^+,e^0,e^-$ as before.

\begin{propos}\label{bjihfebvo}
The  bimodule covariant derivative
$\nabla$ in \ref{kdyrnn} preserves the metric $g$ (as above) if and only if
 $n_\pm$ and $r$ are real, and
$q^{-1}\,g_{++}\,m_+^*
 = g_{00}\,\mu$ and $q^{-1}\,g_{00}\,\nu^*= g_{--}\,m_-$.
\end{propos}
\proof From \ref{ksdzjcv} and the fact that $\extd g_\bullet=0$,
the following matrix must be antiHermitian,
\begin{eqnarray*}
g_\bullet\, \Gamma \,=\, -\,\left(\begin{array}{ccc}g_{++}\,n_+\, e^0 & g_{++}\,m_+\, e^+ & 0 \\
g_{00}\,\mu\, e^- & g_{00}\,r\,e^0 & g_{00}\,\nu\, e^+ \\
0 & g_{--}\,m_-\,e^- & g_{--}\,n_-\, e^0\end{array}\right)\ ,
\end{eqnarray*}
and so,
\begin{eqnarray*}
(g_{++}\,m_+\, e^+)^* &=& -\,q^{-1}\,g_{++}\,m_+^*\,e^-
 \,=\, -\, g_{00}\,\mu\, e^-\ , \cr
(g_{00}\,\nu\, e^+)^* &=&-\,q^{-1}\,g_{00}\,\nu^*\,e^-
\,=\, -\,  g_{--}\,m_-\,e^-\ .\quad\square
\end{eqnarray*}

\begin{propos}\label{cvshk}
The bimodule covariant derivative
$\nabla$ in \ref{kdyrnn}  is metric preserving, torsion compatible and  star compatible if and only if it lies in  the 2-component 1-parameter moduli space given by
 \[ n_+\quad{\rm real},\quad r=\begin{cases}0\\ {2\over q^2-1}\end{cases}\]
and
\[
n_-\,=\,\frac{-\,n_+}{1+(1-q^2)n_+}\ ,\quad m_+\,=\,q^4\,n_+-(1+q^2)\ ,\quad m_-\,=\, {m_+n_-\over q^4n_+},\]
\[\nu= q\,g_{--}\,m_-/g_{00}\ ,\quad \mu\, = \,  q^{-1}\,g_{++}\,m_+/g_{00}\ .\]
The torsion is given by
\begin{eqnarray*}
T(e^\pm)=0\ ,\quad
T(e^0) \,=\,q\,(g_{--}\,m_--g_{++}\,m_+-q^2\,g_{00})\, e^+\wedge e^-/g_{00}\ .
\end{eqnarray*}
\end{propos}
\proof Combine the previous results.  Note that the two choices of $r$ correspond to $\rho(r)=\pm 1$.
\quad$\square$

\begin{theorem} \label{nceicnpncuu} The moduli space of metric-preserving star compatible connections in Proposition~\ref{cvshk} has precisely four points where the connection is torsion free. Of these just one, with $r=0$ and a particular choice of $n_+$,  has a classical limit as $q\to 1$, where it becomes the classical Levi-Civita connection.
\end{theorem}
\proof In Proposition~\ref{cvshk}, for each choice of $r$, the torsion vanishing imposes a quadratic equation on $n_+$, i.e.  the torsion vanishes at either two, one or no points of our moduli space, depending on the discriminant. We shall show that
\begin{eqnarray*}
disc=(g_{00}\, q^2(1-q^2))^2 +(g_{++}+g_{--})^2 + 2\,g_{00}\,q^2\,(1-q^2)\,(g_{--}-g_{++})>0
\end{eqnarray*}
For brevity we write $x=g_{++}$, $y=g_{00}q^2(1-q^2)$ and $z=g_{--}$. Then we want to show
\[ disc=y^2+(x+z)^2+2y(z-x)=(x+y+z)^2-4xy>0.\]
Note that the diagonal metric entries are strictly positive so $x,z>0$. If $q^2>1$ we have $y<0$ so the we are done as $(x+y+z)^2\ge 0$ and $4xy<0$. If $q^2=1$ then $y=0$ and we are also done as $(x+z)^2>0$. It remains to consider $q^2<1$, so $y>0$. Hence the first factor of
\[ disc=(x+y+z+2\sqrt{xy})(x+y+z-2\sqrt{xy})\]
is also strictly positive, while also
\[ x+y+z-2\sqrt{xy}= (\sqrt{x}-\sqrt{y})^2+z >0.\]

Hence $n_+$ has exactly two solutions, given in terms of the metric by
\begin{eqnarray} \label{bcdisuobvo}
n_+ (1-q^2)\,=-1- \frac{g_{--}\,  - g_{++}\,  +g_{00}\, q^2(1-q^2)   \pm\,\sqrt{disc} }
{2\, q^4g_{++}}\ .
\end{eqnarray}
This and our two choices for $r$ gives the four connections. Clearly for a classical limit we need $r=0$ as the other choice is not defined as $q\to 1$. For $n_+$, suppose that $1-q^2$ is small. Then we can choose the square root of the discriminant
to be, to first order in $1-q^2$,
\begin{eqnarray*}
\sqrt{disc}\,\approx\,
(g_{++}+g_{--})\left( 1+\frac{q^2\,g_{00}\, (q^2-1)\,(g_{++}-g_{--})}{(g_{++}+g_{--})^2}
\right).
\end{eqnarray*}
Substituting this into the formula (\ref{bcdisuobvo}) for $n_+$, we see that the root with the $+$
sign gives no limit for $n_+$ as $q\to 1$, but the root $-$ sign has a limit  as $q\to 1$.  This gives our unique conneciton. Its limit corresponds to
\begin{eqnarray*}
n_+ \to 2 - {g_{00}\over g_{++}+g_{--}}
\end{eqnarray*} which  allows us (see Section~6.3) to identify the limit as the classical Levi-Civita connection of the metric when this applies. \quad$\square$

Note that we have been working with Hermitian metrics and in order to speak of the usual Riemannian metric and its Levi-Civita connection we need to consider whether the metric has a real form over $\R$. In our framework this is deferred and we work directly with the Hermitian metric. However, in the classical limit we look at the issue carefully for the metric on $SU_2$ and find that our class of metrics is an actual real metric precisely when $g_{++}=g_{--}$ and the unique connection identified in the theorem above is then its usual Levi-Civita one.

\begin{cor} \label{kzdvvcdijsbc}
There are precisely two points in the moduli space of metric-preserving torsion compatible and star compatible connections in Proposition~\ref{cvshk}  for which $\sigma$ is a braiding, namely
\begin{eqnarray*}
r=\begin{cases}0\\ 2\over q^2-1\end{cases},\quad n_+\,=\,\frac{1+q^2}{q^4}\ ,\quad  n_-\,=\,-(1+q^2)  \ ,\quad  m_+\,=\,m_-\,=\,\nu=\mu=0\ .
\end{eqnarray*}
Then  the eigenvector table for $\sigma$ becomes
\begin{eqnarray*}
q^2\,e_+\tens e_-+e_-\tens e_+   &\mathrm{eigenvalue}&  1\ ,\cr
q^2\,e_+\tens e_--e_-\tens e_+   &\mathrm{eigenvalue}&  -1\ ,\cr
e_0\tens e_0   &\mathrm{eigenvalue}&  \pm 1\ ,\cr
q^2\,e_+\tens e_0+q^{-2}\,e_0\tens e_+   &\mathrm{eigenvalue}&
1\ ,\cr
q^2\,e_+\tens e_0-q^{-2}\,e_0\tens e_+   &\mathrm{eigenvalue}&
-1\ ,\cr
q^4\, e_0\tens e_-+e_-\tens e_0  &\mathrm{eigenvalue}&
1\ ,\cr
q^4\, e_0\tens e_--e_-\tens e_0   &\mathrm{eigenvalue}&
-1\ ,\cr
e_+\tens e_+  &\mathrm{eigenvalue}&  q^2\ ,\cr
e_-\tens e_-  &\mathrm{eigenvalue}&  q^{-2}\ ,
\end{eqnarray*}
with $+1$ and $-1$ respectively for the two choices for $r$.  Such connections can never be torsion free (for a non-degenerate Riemannian structure). The case $r=0$ is the unique such connection having a  limit as $q\to 1$.
\end{cor}
\proof We combine the previous results.\quad$\square$

This result is parallel to Corollary~\ref{corbras3} for the permutation group but the case $r=0$ with classical limit is not the left Maurer-Cartan connection $\nabla^L=0$. However, it is closely related to (but not equal to) the right Maurer-Cartain connection as follows. Recall from Section~4.3 case 2) that for a bicovariant calculus the right Maurer-Cartan form is given in terms of $\nabla^L$ as stated there. Our calculus is {\em not} bicovariant but we do have a right coaction $\rho=(\id\tens\pi)\Delta$ of $\C\Z$ which we have already used above. We also have a $q$-monopole connection form $\varpi_{mon}:\C\Z\to \Omega^1\C_q[SU_2]$ on $\C_q[SU_2]$ as a quantum principal bundle, see \cite{Ma:rsph} for a recent introduction with the 3D calculus. Thus we can write
\[ \nabla^L\xi= -\varpi_{mon}(S^{-1}\xi_{[1]})\tens\xi_{[0]}\]
where now $\rho(\xi)=\xi_{[0]}\tens\xi_{[1]}$.  A short computation shows that this $\nabla^L$ is exactly the $r=0$ choice in Corollary~\ref{kzdvvcdijsbc}. We obtain it as the unique metric-preserving star-compatible connection on $\C_q[SU_2]$ with braiding and a classical limit. The other choice is `purely quantum'  as it does not have a  classical limit as $q\to 1$.

We see that the moduli space in Proposition~\ref{cvshk} of metric preserving torsion compatible and star compatible connections includes a unique $q$-deformed Levi-Civita, and a unique $q$-deformed other connection with torsion, characterised by $\sigma$ obeying the braid relations, and several purely quantum candidates. The connection with torsion does not actually depend on the metric in our class of metrics and is related to the monopole. To complete the picture we also have:

\begin{cor}\label{starpressu2}
There are precisely two points in the moduli space of metric-preserving star compatible connections in Proposition~\ref{cvshk}  for which the connection is star-preserving, namely \[
r=\begin{cases}0\\ 2\over q^2-1\end{cases},\quad n_+\,=\,(\frac{g_{--}}{q^4\,g_{++}}-1)/(1-q^2),\quad n_-=-q^4n_+{g_{++}\over g_{--}}\]
\[ m_+={{g_{--}\over g_{++}}-1\over 1-q^2},\quad m_-=-m_+{g_{++}\over g_{--}},\quad \nu=-q m_+{g_{++}\over g_{00}},\quad \mu=q^{-1}m_+{g_{++}\over g_{00}}\]
Neither of these have a  limit as $q\to 1$ unless $g_{++}=g_{--}$. In this case there is a unique choice $r=0$ with a limit and it coincides with the connection in Corollary~\ref{kzdvvcdijsbc}.
\end{cor}
\noindent {\bf Proof:}\quad We combine Propositions~\ref{bchdsucvasu} and \ref{cvshk}, i.e. we impose the additional condition $\nu=-\mu^*q^2$ in the context of the latter. \quad$\square$

We are thinking of $g_{++}$ and $g_{--}$ as constants independent of $q$, otherwise we only need $g_{++}\to g_{--}$ suitably in the limit. We are led to the same restriction as corresponding at least when $q=1$ to a real metric and we see that in this case our $q$-monopole-related connection has both a braiding (obeying the braid relations) and is fully star-preserving and metric preserving, but has torsion even when $q=1$. Indeed, the  calculations in this section are made under the assumption $q^2-1\ne 0$ and we find a much grater rigidity in this  $q$-deformed world than visible classically at $q=1$ (where braiding and star-preserving are possible with or without torsion).

\subsection{Levi-Civita connections and a reality condition for Riemannian metrics}

The above example allows us to illustrate clearly the difference between the `Hermitian' Riemannian metric we have discussed obeying the symmetry appropriate for a Hilbert $C^*$ module, as in \ref{hermdeff}, and conventional Riemannian geometry. Clearly there is another condition required for a
`real' Riemannian metric: the inner product of two real vectors is real.
In terms of the complex bundles we find ourselves working with, we need some requirement
that the Riemannian metric preserves the star operation on the bundle.

Thus, we look at how in our approach  the classical (i.e.\ $q=1$) Levi Civita connection for the given class of  Riemannian structures on $SU_2$  emerges. We shall assume, as before, that the connection is left
$SU_2$ invariant and right circle invariant. Thus the Christoffel symbols are still of the form
given in \ref{kdyrnn}.

 First we have the strangest condition on the covariant derivative, at least from the point of real classical differential geometry, where it is automatic in terms of real coordinates and real Christoffel symbols; compatibility with the star operation.
 Given that $\sigma$ is just transposition, the
 condition (\ref{bcihdbbdiuv}) for $\nabla$ preserving the star operation can be written in the form
 \begin{eqnarray*}
\nabla(\xi^*)\,=\, \nabla(\xi)^{*\tens *}\ .
\end{eqnarray*}
This gives
\begin{eqnarray*}
\nabla(e^{i*})\,=\, -\Gamma^{i*}_k\tens e^{k*}
\end{eqnarray*}
If we define $\varphi(+)=-$, $\varphi(-)=+$ and $\varphi(0)=0$, then we have (again remembering that
we are in the $q=1$ case)
$e^{i*}=-e^{\varphi(i)}$. Then the condition for $\nabla$ preserving the star operation becomes
\begin{eqnarray*}
-\nabla(e^{\varphi(i)})\,=\, \Gamma^{i*}_k\tens e^{\varphi(k)}\,=\,  \Gamma^{\varphi(i)}_k\tens e^{k}
\end{eqnarray*}
which on simplification gives
\begin{eqnarray*}
\Gamma^{i*}_k\,=\,  \Gamma^{\varphi(i)}_{\varphi(k)}\ .
\end{eqnarray*}
Referring back to the Christoffel symbols in \ref{kdyrnn}, from this equation we can read off
the consequences that $r$ is imaginary, $n^*_+=-n_-$, $m^*_+=-m_-$ and $\mu^*=-\nu$.

Secondly, the metric preserving condition gives the same result as in
\ref{bjihfebvo}, which is that  $n_\pm$ and $r$ are real, and
$g_{++}\,m_+^*
 = g_{00}\,\mu$ and $g_{00}\,\nu^*= g_{--}\,m_-$.
 Combining this with $\nabla$ preserving the star operation, we get $r=0$, $n_+=-n_-$, and
\begin{eqnarray*}
g_{++}\, m_-\,=\,-\,g_{00}\,\mu\,=\, g_{--}\, m_-\ .
\end{eqnarray*}
Consequently (given the metric and star conditions),
 we see if $m_-\neq 0$ that necessarily $g_{++}=g_{--}$.
 Similarly, if $m_-= 0$ we necessarily have $\mu=\nu=0$.

 Thirdly, the torsion is given, as in
\ref{cgvvcsuvc}, by
\[
T(e^\pm)= (m_\pm - n_\pm  \pm 2  )\, e^\pm\wedge e^0
\ ,\quad T(e^0) = (\nu-\mu-1)\,e^+\wedge e^-\ .
\]
If the torsion vanishes, from $T(e^0)=0$
 we cannot have $\mu=\nu=0$, and following the comment above
(given the metric and star conditions) we must have $m_-\neq 0$. Similarly from $T(e^\pm)=0$
we must have $m_\pm$ real. Then $\mu$ and $\nu$ are also real, and from the equations above we find
\begin{eqnarray}\label{bcusuyavc}
\nu\,=\,-\,\mu\,=\,\frac12\ ,\quad m_-\,=\,-\,m_+\,=\, \frac{g_{00}}{2\,g_{++}}\ ,\quad n_+\,=\, -\,n_-\,=\,2-\frac{g_{00}}{2\,g_{++}}\ .
\end{eqnarray}
This gives the classical Levi-Civita connection, in the case that $g_{++}=g_{--}$, and if $g_{++}\neq g_{--}$ we have no star preserving metric preserving torsion-free connection.  This condition is also what we would get at $q=1$ with $e^{\pm}{}^*=-e^\mp$ and requiring our form of $g$ to be real in a hermitian basis.

\subsection{Riemannian geometry for a 2D calculus on a quantum sphere}

We can construct a quantum sphere from $\C_q[SU_2]$ by taking the invariant
elements under the right $\mathbb{CZ}$ coaction. We then take the 3D calculus described earlier, and take the elements which are both horizontal (in this case, an algebra multiple of $e^+$
or $e^-$) and invariant to the
$\mathbb{CZ}$ coaction. This we take to be the differential calculus on the
quantum sphere, as part of the standard construction for the $q$-monopole cf\cite{MaBrz:gau}; see \cite{Ma:rsph, ncsheaf} for recent works. In this construction the horizontal 1-forms coincide with the kernel of the
`vertical vector fields' map
\begin{eqnarray*}
\Omega^1 \C_q[SU_2] \to \C_q[SU_2] \tens \Omega^1 \mathbb{CZ}\ .
\end{eqnarray*}
given by differentiating the $\mathbb{CZ}$ coaction. Now we look for connections on $\C_q[SU_2]$ which descend to connections on the sphere:

We consider which bimodule covariant derivative
$\nabla$ in Proposition~\ref{kdyrnn}  on $\C_q[SU_2]$  descend to the quantum  sphere.
This requires that horizontal 1-forms are sent to elements of $\Omega^1$ tensor
horizontal 1-forms by $\nabla$. Combining this with the bimodule condition and
right $\mathbb{CZ}$ invariance gives
\begin{eqnarray*}
\nabla^L(e^\pm)\,=\, n_\pm\, e^0\tens e^\pm\ .
\end{eqnarray*}
(The value of $\nabla^L(e^0)$ is not relevant.)
The values of $n_\pm$ are calculated so that $\nabla(x\, e^+)$
and $\nabla(y\, e^-)$ are horizontal when $x\, e^+$ and $y\, e^-$ are
invariant forms under the right $\mathbb{CZ}$ coaction.
 We do this by a bit of brute force, considering
all the possibilities (up to a right invariant element of the algebra):
\begin{eqnarray*}
\extd(ac) \,=\, (1+q^2)\, ac.e^0+\dots &,& \extd(a^2) \,=\, (1+q^2)\, a^2.e^0+\dots \ ,\cr
\extd(c^2) \,=\, (1+q^2)\, c^2.e^0+\dots &,& \extd(d^2) \,=\, -q^{-4}(1+q^2)\, d^2.e^0
+\dots \ ,\cr
\extd(b^2) \,=\, -q^{-4}(1+q^2)\, b^2.e^0+\dots &,&
 \extd(bd) \,=\, -q^{-4}(1+q^2)\, bd.e^0+\dots \ .
\end{eqnarray*}
Here the dots indicate a horizontal form (a multiple of $e^+$
or $e^-$).
It follows that
\begin{eqnarray*}
n_-\,=\, -(1+q^2)\ ,\quad n_+\,=\, q^{-4}(1+q^2)\ .
\end{eqnarray*}
This gives the values in \ref{kzdvvcdijsbc}. On the sphere the generator $e^0$ does not appear,
so the connection is torsion free and unique among covariant derivatives that can be obtained in this way. It preserves the Hermitian metric given by
restricting the Hermitian metric to $\C_q[SU_2]$. This is also the same covariant derivative
as the one given by the second author in \cite{Ma:rief}, as follows easily by comparison with the explicit formula for the action of the covariant derivative $D$ on sections in the proof of Theorem~5.1 there.

\end{document}